\documentclass[12pt]{amsart}
\usepackage{amsfonts}
\usepackage{amsmath,amssymb,amsthm,latexsym,graphicx,graphics}
\usepackage{multicol}
\usepackage{color}
\usepackage{url,hyperref}
\usepackage{euscript}
\usepackage{wrapfig}
\usepackage{caption}
\usepackage{url,hyperref}
\usepackage{euscript,color}

\newtheorem{teo}{Theorem}

\newtheorem*{lemn}{Lemma}
\newtheorem{prp}{Proposition}

\newtheorem{cor}{Corollary}

\newtheorem*{proofn}{Proof}

\begin{document}

\title[M\"{o}bius fluid dynamics on the unitary groups]{M\"{o}bius fluid dynamics on the unitary groups}
\author{Daniela Emmanuele$^{1}$, Marcos Salvai$^{2}$ and Francisco Vittone$^{3}$}
%
%
\begin{abstract}
We study the nonrigid dynamics induced by the standard birational actions of
the split unitary groups $G=O_{o}\left( n,n\right) $, $SU\left( n,n\right) $
and $Sp\left( n,n\right) $ on the compact classical Lie groups $M=SO_{n}$, $%
U_{n}$ and $Sp_{n}$, respectively. More precisely, we study the geometry of $%
G$ endowed with the kinetic energy metric associated with the action of $G$
on $M,$ assuming that $M$ carries its canonical bi-invariant Riemannian
metric and has initially a homogeneous distribution of mass. By the least
action principle, force free motions (thought of as curves in $G$)
correspond to geodesics of $G$. The geodesic equation may be understood as
an inviscid Burgers equation with M\"{o}bius constraints. We prove that the
kinetic energy metric on $G$ is not complete and in particular not
invariant, find symmetries and totally geodesic submanifolds of $G$ and
address the question under which conditions geodesics of rigid motions are
geodesics of $G$. Besides, we study equivalences with the dynamics of
conformal and projective motions of the sphere in low dimensions.

\end{abstract}
\maketitle


\section{Introduction and statements of the results}
\let\thefootnote\relax\footnotetext{$^1$ emman@fceia.unr.edu.ar; Universidad Nacional de Rosario, Pellegrini 250, 2000, Rosario, Argentina. }
\let\thefootnote\relax\footnotetext{$^{2}$ salvai@famaf.unc.edu.ar; FaMAF, Universidad Nacional de C\'ordoba - CIEM, CONICET, 
	Ciudad Universitaria, (5000) C\'{o}rdoba, Argentina. }
\let\thefootnote\relax\footnotetext{$^{3}$ vittone@fceia.unr.edu.ar; Universidad Nacional de Rosario - CONICET, Pellegrini 250, 2000, Rosario, Argentina. }
Geometric mechanics, introduced by Arnold in \cite{ArnoldArt}, describes the
force free motions of a system as geodesics of a Riemannian metric on the
configuration space, given by the kinetic energy. Among the first examples
we find invariant metrics on $SO_{3}$ and on the diffeomorphism group Diff$%
_{\mu }\left( M\right) $ of a compact Riemannian manifold $M$ preserving the
volume form $\mu $, in order to describe, respectively, force free motions
of the rigid body and the dynamics of an ideal flow on $M$.

This approach had an even greater impact when many classical differential
equations were recognized within this theory --that is, as geodesic
equations--, among others, Korteweg-de Vries \cite{OK87} (see more examples
for instance in \cite{Viz08, BBM14}).

In this paper we study the dynamics of certain finite dimensional Lie groups
endowed with the kinetic energy metric induced by the action of them on
Riemannian manifolds, which do not preserve the volume. We begin by putting
the situation in a broader context. 
Ebin and Marsden \cite{EbinMarsden} extended Arnold's ideas for the whole
diffeomorphism group Diff$\left( M\right) $, studying in this way
compressible flows. In this case, the metric is no longer invariant. The
force free motions $t\mapsto\gamma\left( t\right) $, or equivalently, the
geodesics in Diff$\left( M\right) $, are rendered in a simple manner:\ $%
t\mapsto\gamma\left( t\right) \left( q\right) $ turn out to be geodesics in $%
M$ for all $q\in M$. This description is no longer valid if the flow is
constrained, for instance in the case of Diff$_{\mu}\left( M\right) $,
which, with the metric inherited from that of Diff$\left( M\right) $, is not
totally geodesic.

Geodesics of the induced metric have also been considered for groups of
diffeomorphisms preserving other structures on $M$, for example for the
group of symplectomorphisms \cite{EbinSympl, Kh2} and contactomorphisms \cite%
{EbinContact}. For these groups, and in general for subgroups of Diff$%
_{\mu}\left( M\right) $, the kinetic energy metric is invariant and has been
studied extensively. In the finite dimensional case, invariant metrics
appear in the rigid context, and force free motions, that is, geodesics,
have been much explored and relatively well understood.

Concerning non-invariant metrics, they arise for compressible flows (the
basic example being Diff$\left( M\right) $, which we mentioned above) and
have been less studied. Bauer and Modin have recovered the interest in this
subject in \cite{BauerModinSemi}.

By non-rigid dynamics we understand the study of a dynamical system given by
a kinetic energy metric on a Lie group, induced by the action of it on a
Riemannian manifold, which is not by isometries. Among the reasons why the
nonrigid dynamics in finite dimension has been addressed to a lesser extent,
we can mention that the most significant actions of noncompact Lie groups on
compact Riemannian manifolds are in general scarce, for instance, for the
important family of compact symmetric spaces: the celebrated Nagano's
Theorem \cite{Nagano} states that if a Lie group $G$ acts on a compact
symmetric space $M$ and contains the identity connected component of its
isometry group, then $M$ is essentially a symmetric $R$-space (or
equivalently, it is extrinsically symmetric) and $G$ is the so called big
group, or group of basic transformations in Takeuchi's terminology in \cite%
{Ta}. The best known examples are the conformal and projective
transformations of the sphere $S^{n}$, the biholomorphic maps of the complex
quadrics and the diffeomorphisms on Grassmannians which are induced by
linear isomorphisms. More generally, the $R$-spaces, homogeneous models of
the important parabolic geometries \cite{CapSlovak, Sharpe}, admit actions
of finite dimensional Lie groups containing the isometries.

The nonrigid (and non-invariant) dynamics of the simplest among the
symmetric $R$-spaces, the sphere, is studied in \cite{Sal} and \cite{LSW},
for its two big groups, associated with conformal and projective
transformations. Also in the nonrigid finite dimensional setting, we can
cite \cite{Burov} (force free affine motions in Euclidean space) and the
recent article \cite{Sideris}, where a non-invariant metric on $SL\left(3,%
\mathbb{R}\right)$ appears in relation to affine motions of 3D ideal fluids
surrounded by vacuum.

A further fundamental problem in this area, apart from existence and
regularity of the geodesic flow, is its completeness, related with the onset
of shocks. The first two issues (existence and regularity) are of interest
only in infinite dimension. The third one (completeness), also in finite
dimension (of course, only when the configuration space is not Riemannian
homogeneous), and arises naturally, as we have pointed out, when the group
acts without preserving the volume.

In finite dimension, in the invariant case, the kinetic energy metric (as
any invariant metric) is complete. Since for infinite dimension the
Hopf-Rinow theorem is no longer valid, not even for invariant metrics,
completeness depends on regularity. Bruveris and Vialard deal with this
topic in depth in \cite{BruverisVialardJEMS}. In \cite{BauerModinSemi}
sufficient conditions for completeness are found in the semi-invariant case.

In this article we study the dynamics of M\"obius motions of the compact
classical Lie groups, mainly drawing on the symmetries of the systems, in
the interest of advancing the understanding of non-rigid motions in finite
dimension. In this context, non-invariant metrics on Lie groups arise
naturally.

We consider the unitary groups $SO_{n}$, $U_{n}$ and $Sp_{n}$, over the
real, complex and quaternionic numbers (that is, the identity components of
the automorphism groups of the positive definite Hermitian forms over these
scalar fields, see Subsection \ref{unitary}).

We recall the standard birational action $\ast $ of the split unitary group $%
G=O_{o}\left( n,n\right) $, $SU\left( n,n\right) $ or $Sp\left( n,n\right) $
on the compact classical Lie group $M=SO_{n}$, $U_{n}$ or $Sp_{n}$,
respectively:%
\begin{equation}
	\ast :G\times M\rightarrow M\text{, \ \ \ \ \ \ \ }\left( 
	\begin{array}{cc}
		A & B \\ 
		C & D%
	\end{array}%
	\right) \ast U=(AU+B)(CU+D)^{-1}\text{,}  \label{mainAction}
\end{equation}%
which we call the \textbf{M\"{o}bius action }of $G$ on $M$. This is one of
the actions studied in \cite{Nagano}, which we mentioned above (see also 
\cite{Th}). Geometrically, it comes from the identification of $M$ with the
Grassmannian of maximal isotropic subspaces, assigning to each operator $U$
in $M$ its graph from the second to the first factor (see Subsection \ref%
{identgrass} below) and the natural action of $G$ on this Grassmannian:\ if%
\begin{equation*}
	\left( 
	\begin{array}{cc}
		A & B \\ 
		C & D%
	\end{array}%
	\right) \left( 
	\begin{array}{c}
		Ux \\ 
		x%
	\end{array}%
	\right) =\left( 
	\begin{array}{c}
		AUx+B \\ 
		CUx+D%
	\end{array}%
	\right) =\left( 
	\begin{array}{c}
		U^{\prime }y \\ 
		y%
	\end{array}%
	\right)
\end{equation*}%
for $x,y=y\left( x\right) \in \mathbb{F}^{n}$, then $U^{\prime }$ is the
right hand side of (\ref{mainAction}).

Up to coverings and connected components, $G$ is the big group of $M$, which
is characterized in \cite{STohoku} as the group of diffeomorphisms of $M$
sending circles in circles (certain distinguished curves in $M$ determined
by three generic points).

Now, we consider the Riemannian structure on $M$ induced by the inclusion in
the space $\mathbb{F}^{n\times n}$ of $\left( n\times n\right) $-matrices
with entries in $\mathbb{F}={\mathbb{R}},\,{\mathbb{C}}$ or the quaternions $%
{\mathbb{H}}$, endowed with the real inner product given by $\left\langle
X,Y\right\rangle =\operatorname{Re}\left(\operatorname{tr}\left( \overline{X}^{T}Y\right)
\right)$, where the superscript $T$ denotes transposition (beware of the
peculiar relationship of the multiplication of quaternionic matrices with
transposition and conjugation; see Theorem 4.1 in \cite{Zh}). It turns out
that the metric on $M$ is bi-invariant and also that M\"{o}bius maps are not
necessarily isometries.

We are concerned with the nonrigid dynamics induced by this non-isometric
action. We take $G$ as the \textbf{configuration space}. This means that, if 
$M$ has initially a homogeneous mass distribution of constant density $1$,
the points of $M$ are allowed to move only in such a way that two
configurations differ in an element of $G$. We will see that the action of $G
$ on $M$ is almost effective (only a finite number of elements act as the
identity), except for $G=O_{o}\left( 2,2\right) $. For this reason, unless
otherwise stated, from now on we assume that $n>2$ when $M=SO_{n}$.
Nevertheless, for the sake of completeness, we describe the M\"{o}bius
action of $O_{o}\left( 2,2\right) $ on $SO_{2}$ in Proposition \ref{O22}.

A curve $\gamma $ in $G$ may be thought of as a \textbf{motion} of $M$. The
Riemannian metric on $M$ allows us to consider the kinetic energy of a
motion and this induces a Riemannian metric on $G$, introduced originally by
V.I. Arnold in \cite{ArnoldArt}, called the \textbf{kinetic energy metric}
on $G$ associated with the action on $M$, with the property that \textbf{%
	force free motions} (thought of as curves in $G$) correspond to \textbf{%
	geodesics} of $G$, via the least action principle.

For each $g\in G$ and $X\in T_{g}G$, the square norm of $X$ for the kinetic
energy metric turns out to be%
\begin{equation*}
	\left\Vert X\right\Vert ^{2}=\int_{M}\left\vert \left. \frac{d}{dt}%
	\right\vert _{0}\left( ge^{tX}\ast q\right) \right\vert ^{2}d\mu \left(
	q\right) \text{,}
\end{equation*}%
where $\left\vert v\right\vert $ denotes the norm of a tangent vector to $M$
(see Subsection \ref{seckinetic}).

The kinetic energy metric on $G$ is the restriction to $G$ of the usual weak
Riemannian metric on $\operatorname{Diff}\left(M\right) $ and the equation for
geodesics may be understood as the \textbf{inviscid Burgers equation with M%
	\"{o}bius constraints} (see \S\ 4.1 in \cite{Khesin}).

Next we state our first result on the geometry of $G$, its incompleteness.
We denote by $I_{n}$ and $0_{n}$ the identity and the zero $\left( n\times
n\right) $-matrices, respectively.

\begin{teo}
\label{geodfinita copy(1)}Let $G$ be endowed with the kinetic energy metric.
Then the curve $\gamma :\mathbb{R}\rightarrow G$ defined by
\begin{equation*}
\gamma (t)=\exp \left( t\left(
\begin{array}{cc}
0_n & I_{n} \\
I_{n} & 0_n%
\end{array}%
\right) \right) =\left(
\begin{array}{cc}
\cosh t~I_{n} & \sinh t~I_{n} \\
\sinh t~I_{n} & \cosh t~I_{n}%
\end{array}%
\right)
\end{equation*}%
is the reparametrization of an inextendible geodesic in $G$ of finite
length. In particular, $G$ is not complete and so the kinetic energy metric
on $G$ is neither left nor right invariant.
\end{teo}

Incompleteness, in a classical context, is related to the existence of
shocks (collapses in finite time).

We have mentioned above the work of Bruveris and Vialard
\cite{BruverisVialardJEMS} about completeness in the invariant case. For the semi-invariant case, in the recent paper \cite{BauerModinSemi},
Bauer and Modin give \emph{sufficient} conditions for completeness. Notice
that the metric on $G$ in this article
is semi-invariant in the sense of \cite{BauerModinSemi}, as we will see in Proposition \ref{CompactoMaximal}.

The following proposition shows the behavior of the mass when $%
\lim_{t\rightarrow \pm \infty }\gamma \left( t\right) $.

\begin{prp}
\label{acum}If $M\neq SO_{n}$ with odd $n$, under the action of $\gamma $,
the mass in $M$ concentrates on $\pm I_{n}$ as $t\rightarrow \pm \infty $.
More precisely, for $\varepsilon =\pm 1$,
\begin{equation*}
\lim_{t\rightarrow \varepsilon \infty }\gamma \left( t\right) \ast
q=\varepsilon I_{n}
\end{equation*}%
for all $q\in \left\{ p\in M\mid p\text{ has no eigenvalue }-\varepsilon
\right\} $, which is an open and dense set in $M$. If $M=SO_{n}$ with odd $n$, then the same property holds for $\varepsilon =1$.
\end{prp}

Next we present some totally geodesic submanifols $N$ of $G$ endowed with
the kinetic energy metric. In our dynamical setting, this means that if a
force free motion begins tangent to $N$, then it remains in $N$.

\begin{teo}
\label{TotGeod}The following inclusions determine totally geodesic
submanifolds \emph{(}in any case the ambient group carries the kinetic
energy metric and subgroup, the induced metric\emph{):}
\begin{equation*}
\begin{tabular}{lll}
\emph{a)} $O_{o}\left( n,n\right) \subset SU\left( n,n\right) $, & \emph{b)}
$U\left( n,n\right) \subset Sp\left( n,n\right) $, & \emph{c)} $O_{o}\left(
n,n\right) \subset Sp\left( n,n\right) $, \\
\emph{d)} $U(n,n)\subset O_{o}(2n,2n)$, & \emph{e)} $Sp\left( n,n\right)
\subset O_{o}(4n,4n)$. &
\end{tabular}%
\end{equation*}
\end{teo}

Totally geodesic subgroups when the kinetic energy metric is invariant are
addressed in \cite{Viz02, ModinPerl}.

\smallskip

Now we consider the maximal compact subgroups $K=SO_{n}\times SO_{n}$, $%
S(U_{n}\times U_{n})$ and $Sp_{n}\times Sp_{n}$ of $G$ (for the three
different possibilities for $G$, namely $O_{o}\left( n,n\right) $, $SU\left(
n,n\right) $ and $Sp\left( n,n\right) $). Up to coverings and connected
components, $K$ is the isometry group of $M$, or, in our dynamical context,
the group of rigid transformation of $M$. 

We comment, without trying to be precise, that $M$ appears, on the one hand,
in a passive manner, as the manifold containing the fluid whose motion is
subject to M\"{o}bius constraints, and also, on the other hand, in an active
way, acting on itself on the right and on the left. 
Notice that $K$ is equal to $%
M\times M$ or is contained in it and that the M\"{o}bius action (\ref%
{mainAction}) restricted to $K$ coincides essentially with the canonical
action $\operatorname{diag}\left( A,D\right) \ast U=AUD^{-1}$ of $M\times M$ on $M$
by isometries (here $\operatorname{diag}(A,D)$ denotes the matrix with diagonal
blocks $A$ and $D$).

Although the kinetic energy metric on $G$ is not invariant, it has a certain
symmetry given by the action of $K\times K$, as the following proposition
states, placing the topic in the framework of \cite{BauerModinSemi}.

\begin{prp}
\label{CompactoMaximal}\emph{a)} The action of $K\times K$ on $G$ given by $%
\left( k_{1},k_{2}\right) \left( g\right) =k_{1}gk_{2}^{-1}$ is by
isometries of the kinetic energy metric of $G$.

\smallskip

\emph{b)} The Riemannian metric on $K$ induced by the kinetic energy metric
of $G$ is bi-invariant; moreover, it is the restriction to $K$ of a product
metric on $M\times M$ such that the map $\rho :K\rightarrow K$ given by $%
\rho \left( g,h\right) =\left( h,g\right) $ is an isometry.
\end{prp}

\smallskip

Since the metric on $K$ is bi-invariant, the geodesics $\sigma $ in $K$
through the identity have the form $\sigma \left( t\right) =\exp \left(
tX\right) $ with $X\in \mathfrak{k}$, where $\mathfrak{k}$ denotes the Lie
algebra of $K$. The following theorem answers partially the question whether
they are also geodesics in $G$. We denote by $\mathfrak{m}$ the Lie algebra
of $M$.

\begin{teo}
\label{XceroceroY}Let $G$ be endowed with the kinetic energy metric and
suppose that $Z\in \mathfrak{k}$ equals $\operatorname{diag}(X,0_n)$ or $\operatorname{diag}%
(0_n,X)$, where $X\in \mathfrak{m}$. Then the geodesic $\gamma (t)=\exp (tZ)$
of $K$ is also a geodesic of $G$.
\end{teo}

\medskip

In \cite{Sal} and \cite{LSW} similar problems as in this article have been
studied for the conformal and the projective maps of spheres. We recall the
actions%
\begin{equation*}
\ast _{p}:SL_{m+1}\left( \mathbb{R}\right) \times S^{m}\rightarrow S^{m}%
\text{,\ \ \ \ \ \ \ \ }\ast _{c}:O_{o}\left( 1,m+1\right) \times
S^{m}\rightarrow S^{m}\text{,}
\end{equation*}%
on $S^{m}$ by projective and conformal transformations, respectively, which
are given by%
\begin{equation}
A\ast _{p}z=Az/\left\vert Az\right\vert \text{ \ \ \ \ \ \ and \ \ \ \ \ \ }%
A\ast _{c}z=z^{\prime }\text{,}  \label{projective}
\end{equation}%
where $z^{\prime }$ is the unique point in $S^{m}$ such that%
\begin{equation}
A\left( 1,z\right) ^{T}\in \mathbb{R}\left( 1,z^{\prime }\right) ^{T}\text{.}
\label{conformal}
\end{equation}

It is well known that in low dimensions there are isomorphisms (up to
coverings) between $SL_{n}\left( \mathbb{R}\right) $ or $O_{o}\left(
1,n\right) $ and some of the groups $G$ that have been considered in this
paper (see for instance \cite{HelG}, pp 517--528). In Section \ref{last} we explore whether the corresponding actions on
spheres or unitary groups are equivalent (the task turned out to be more
demanding than expected). This will have implications in the dynamics of
M\"obius motions of the compact classical Lie groups: Translating the
results of \cite{Sal} and \cite{LSW}, Corollary \ref{coro} provides a
refinement of Theorem \ref{XceroceroY} in some low dimensional cases.

\section{Preliminaries and basic facts}

\subsection{Unitary groups\label{unitary}}


In this section we recall some basic facts about the classical compact Lie
groups. For more details on the subject we refer to \cite{Har} and \cite{Th}%
. Let $\mathbb{F}=\mathbb{R}$, ${\mathbb{C}}$ or ${\mathbb{H}}$ (the skew
field of quaternions) and let $N\in\mathbb{N}$. Then $\mathbb{F}^{N}$ is an $%
\mathbb{F}$-vector space, where in the case of $\mathbb{F}={\mathbb{H}}$, we
choose it to be a right ${\mathbb{H}}$-vector space. As usual, we denote by $%
GL(N,\mathbb{F})$ the group of invertible $\mathbb{F}$-linear endomorphisms
of $\mathbb{F}^{N}$.

For $n\in \mathbb{N}$ and $m\in \mathbb{N}\cup \left\{ 0\right\} $ set $%
b_{n,m}^{\mathbb{F}}:\mathbb{F}^{n+m}\times \mathbb{F}^{n+m}\rightarrow
\mathbb{F}$ the sesquilinear form
\begin{equation}
b_{n,m}^{\mathbb{F}}((x,y),(\xi,\eta))=\overline{x}^{T}\,\xi-\overline{y}%
^{T}\,\eta\text{,}  \label{form}
\end{equation}%
where $x,\,\xi\in \mathbb{F}^{n}$, and $y,\,\eta\in \mathbb{F}^{m}$ are
column vectors, the bar denotes conjugation and the superscript $T$
indicates transposition, as before.

The form $b_{n,m}^{\mathbb{F}}$ is called the $\mathbb{F}$\emph{-Hermitian
symmetric inner product of signature }$(n,m)$ on $\mathbb{F}^{n+m}$ and the
space $\mathbb{F}^{n+m}$ endowed with this form is called $\mathbb{F}^{n,m}$%
. Let $U(n,m,\mathbb{F})$ be the Lie group
\begin{equation*}
U(n,m,\mathbb{F})=\{A\in GL(n+m,\mathbb{F})\mid b_{n,m}^{\mathbb{F}%
}(Ax,Ay)=b_{n,m}^{\mathbb{F}}(x,y)\ \ \text{for all }\,x,y\in \mathbb{F}%
^{n+m}\}\text{.}
\end{equation*}%
We shall denote by $U_{o}(n,m,\mathbb{F})$ the connected component of the
identity of $U(n,m,\mathbb{F})$. For $\mathbb{F}={\mathbb{R}}$ or ${\mathbb{C%
}}$, $SU(n,m,\mathbb{F})=\left\{ A\in U(n,m,\mathbb{F})\mid {\mathrm{det}}
A=1\right\} $.

We are interested in the cases where $m=n$ or $m=0$. For brevity, we call $%
U\left( n,\mathbb{F}\right) =U(n,0,\mathbb{F})$. Then, for each $n\in
\mathbb{N}$ we have:

- the orthogonal, unitary and symplectic groups $U(n,{\mathbb{R}})=O_{n}$, $%
U(n,\mathbb{C})=U_{n}$ and $U(n,\mathbb{H})=Sp_{n}$ (also called $HU\left(
n\right) $),

- $U(n,n,{\mathbb{R}})=O(n,n)$, $U(n,n,{\mathbb{C}})=U(n,n)$ and $U(n,n,{%
\mathbb{H}})=Sp\left( n,n\right) $ (also called $HU\left( n,n\right) $).

\subsection{The action of $U(n,n,\mathbb{F})$ on $U(n,\mathbb{F})$\label%
{identgrass}}


An $\mathbb{F}$-subspace $V$ of $\mathbb{F}^{n,n}$ is said to be \textbf{%
totally isotropic} if $\left. b_{n,n}^{\mathbb{F}}\right\vert _{V\times
V}\equiv 0$. A maximal totally isotropic subspace must
have dimension $n$ (see \cite{Har}). Denote by Grass$_{0}(n,\mathbb{F})$ the Grassmannian of
maximal totally isotropic subspaces of $\mathbb{F}^{n,n}$. Then, it is not
difficult to see that the map
\begin{equation}
f:U(n,\mathbb{F})\rightarrow \text{Grass}_{0}(n,\mathbb{F}),\ \ \ \ \ \ \
f(T)=\operatorname{rgraph}(T)=\{(Tx,x))\mid x\in \mathbb{F}^{n}\}  \label{Ugrasss}
\end{equation}%
is a diffeomorphism, where $\operatorname{rgraph}(T)$ is the graph of $T$ reflected
with respect to the diagonal, or equivalently, the graph of $T^{-1}$, or the
graph of $T$ as a map from the second factor into the first.

Now, there is a canonical action of $U(n,n,\mathbb{F})$ on Grass$_{0}(n,%
\mathbb{F})$. The identification $f$ induces an action of $U(n,n,\mathbb{F})$
on $U\left( n,\mathbb{F}\right) $, as follows: Given $U\in U(n,\mathbb{F})$,
let $V=\operatorname{rgraph}(U)$. So, if $T\in U(n,n,\mathbb{F})$ has the matrix
form
\begin{equation*}
T=\left(
\begin{array}{cc}
A & B \\
C & D%
\end{array}%
\right) \text{,}
\end{equation*}%
then $T(V)=\{(AUx+Bx,CUx+Dx)\mid \,x\in \mathbb{F}^{n}\}$. On the other
hand, there exists $U^{\prime }\in U(n,\mathbb{F})$ such that $T(V)=\operatorname{%
rgraph}(U^{\prime })=\{(U^{\prime }y,y)\mid \,y\in \mathbb{F}^{n}\}$ and so
we must have $y=(CU+D)x$ and consequently the action reads, as in (\ref%
{mainAction}),
\begin{equation}
T\ast U=U^{\prime }=(AU+B)(CU+D)^{-1}\text{.}  \label{action}
\end{equation}%
By Witt's Theorem, the canonical action of $U(n,n,\mathbb{F})$ on Grass$%
_{0}(n,\mathbb{F})$ is transitive, and then so is the action restricted to $%
U_{o}(n,n,\mathbb{F})$. Through $f$ we obtain a transitive action of $%
U_{o}(n,n,\mathbb{F})$ on $U(n,\mathbb{F})$.

At this point we should make some remarks:

- As it is well known, $O(n,n)$ is not connected (see \cite{ONeill}). However, since it acts
transitively on $O_{n}$, then $O_o(n,n)$ acts transitively on $%
(O_{n})_{o}=SO_n$.

- The action of $U(n,n,{\mathbb{C}})=U(n,n)$ on $U_{n}$ is not effective. It
is not difficult to see that actually, for each $u\in S^{1}\subset {\mathbb{C%
}}$ and $T\in U\left( n,n\right) $, we have $(uT)\ast U=T\ast U$. In order
to obtain an effective action, we can take the group $PU(n,n)=U(n,n)/S^{1}$.
For convenience, we will consider the action of the matrix group $SU(n,n)$
on $U_{n}$, which is almost effective since for $T$ in this group, $T\ast
U=U $ for all $U\in U_{n}$ if and only if $\ T=uI_{2n}$ with $u^{2n}=1$.

We comment that the group $G=U\left( n,n,\mathbb{F}\right) $ acts
canonically on the Grassmannian of maximal space-like subspaces of $\mathbb{F%
}^{n,n}$, from which $M=U\left( n,\mathbb{F}\right) $, via the
identification (\ref{Ugrasss}), is a border inheriting the action.

\subsection{The kinetic energy metric\label{seckinetic}}


In this section we review the restriction to a finite dimensional subgroup
of the usual $L^{2}$ weak Riemannian metric on the space of all
diffeomorphisms of a compact Riemannian manifold (see for instance Theorem
9.1 (i) in \cite{EbinMarsden}), setting in a broader context some of the
concepts presented in \cite{Sal} and \cite{LSW}.

Let $G$ be a connected Lie group acting smoothly and almost effectively on a
compact oriented Riemannian manifold $M$. Suppose that $M$ has initially a
homogeneous distribution of mass of constant density $1$ and that the points
are allowed to move only in such a way that two configurations differ in an
element of $G$. Then the configuration space may be naturally identified
with $G$ (aside from the fact that a finite number of elements act by the
identity) and a curve $\gamma $ in $G$ may be thought of as a motion of $M$.

The total kinetic energy $E(\gamma )$ of the motion $\gamma $ at the instant
$t$ is given by
\begin{equation}
E(t)=\frac{1}{2}\int_{M}\rho _{t}(q)|v_{t}(q)|^{2}d\mu,  \label{energy}
\end{equation}%
where $d\mu $ is the volume form on $M$ and if $q=\gamma (t)(p)$ for some $%
p\in M $, then
\begin{equation*}
v_{t}(q):=\left. \frac{d}{ds}\right\vert _{t}\gamma (s)(p)\in T_{q}M,\ \ \ \
\ \ \ \rho _{t}(q)=1/{\mathrm{det}}(d\gamma (t)_{p})
\end{equation*}%
are the velocity of the point $q$ and the density at $q$ at the instant $t$,
respectively. This is the smooth analogue of $\sum_{i}\frac{1}{2}%
m_{i}v_{i}^{2}$, once one recognizes $\rho d\mu $ as the mass element.

Changing variables in (\ref{energy}), one obtains
\begin{equation}
E(t)=\frac{1}{2}\int_{M}\left\vert \left. \frac{d}{ds}\right\vert _{t}\gamma
(s)(p)\right\vert ^{2}d\mu .  \label{energy2}
\end{equation}

Based on the principle of least action, we say that a smooth curve $\gamma $
in $G$, thought of as a motion of $M$, is force free if it is a critical
point of the kinetic energy functional.

In order to introduce a Riemannian metric on $G$ whose geodesics correspond
to force free motions of $G$, for $g\in G$ and $X\in T_{g}G$, we define the
map $\widetilde{X}:M\rightarrow TM$ by
\begin{equation*}
\widetilde{X}(q)=\left. \frac{d}{dt}\right\vert _{0}\gamma (t)(q)\in
T_{g(q)}M,
\end{equation*}%
where $\gamma $ is any smooth curve in $G$ with $\gamma (0)=g$ and $\gamma
^{\prime }(0)=X$. Standard arguments show that the map $\widetilde{X}$ is
well-defined and smooth, and that $\widetilde{X}$ is a vector field on $M$
if and only if $X\in T_{e}G$. Following \cite[Prop. 1]{Sal}, one has that
\begin{equation}
\langle X,Y\rangle =\int_{M}\langle \widetilde{X}(q)\,,\,\widetilde{Y}%
(q)\rangle d\mu \text{,}  \label{dinmet}
\end{equation}%
for $X,Y\in T_{g}G$ and $g\in G$, defines a Riemannian metric on $G$ such
that a curve $\gamma $ in $G$ is a geodesic if and only if it is force free
(on the basis of the principle of least action). This metric is called the
kinetic energy metric of $G$ induced by the action on $M$.

\section{Symmetries of the kinetic energy metric}

From now on, we will consider $G=O_{o}(n,n),\,SU(n,n)$ or $Sp\left(
n,n\right) $, their maximal compact subgroups $K=SO_{n}\times SO_{n}$, $%
S(U_{n}\times U_{n})$ or $Sp_{n}\times Sp_{n}$, and $M=SO_{n},\,U_{n}$ or $%
Sp_{n}$, respectively.

\medskip
\noindent
\textbf{Proof of Proposition \protect\ref{CompactoMaximal}. }
{\textit
a) For $k_{1},k_{2}\in K$, set
\begin{equation}
F_{k_{1},k_{2}}:G\rightarrow G,\ \ \ \ \ g\mapsto
F_{k_{1},k_{2}}(g)=k_{1}gk_{2}{}^{-1}\text{.}  \label{Fk1k2}
\end{equation}
We have to verify that $F_{k_{1},k_{2}}$ is an isometry of the kinetic
energy metric on $G$.}

\textit{
Let $X\in T_{g}G$ and let $\gamma $ be a curve in $G$ such that $\gamma
(0)=g $, $\gamma ^{\prime }(0)=X$. Then $dF_{k_{1},k_{2}}(X)=\left. \frac{d}{%
dt}\right\vert _{0}k_{1}\gamma (t)k_{2}^{-1}$ and so for each $q\in M$ we
have that%
\begin{equation*}
\widetilde{dF_{k_{1},k_{2}}(X)}(q)=\left. \frac{d}{dt}\right\vert
_{0}(k_{1}\gamma (t))\ast (k_{2}^{-1}\ast q)\text{.}
\end{equation*}
}
\textit{
Observe that since $k_{i}\in K$, then $\displaystyle k_{i}=\operatorname{diag}\left(
A_{i},B_{i}\right) $ with $A_{i},B_{i}\in M$, and so $k_{i}\ast
q=A_{i}qB_{i}^{-1}=R_{B_{i}^{-1}}L_{A_{i}}(q)$. Define
\begin{equation*}
f:M\rightarrow \mathbb{R}\text{,\ \ \ \ \ \ \ }f(q)=\left\vert \left. \frac{d%
}{dt}\right\vert _{0}(k_{1}\gamma (t))\ast q\right\vert ^{2}\text{,}
\end{equation*}
and $h:M\rightarrow M$ by $h=R_{B_{2}}L_{A_{2}^{-1}}$. Then,
\begin{equation}
\left\Vert dF_{k_{1},k_{2}}X\right\Vert ^{2}=\int_{M}\left(f\circ
h\right)\,d\mu =\int_{M}f(h^{\ast }d\mu )=\int_{M}fd\mu,  \label{nodF}
\end{equation}%
since $d\mu $ is bi-invariant by the action of $M$.}

\textit{
On the other hand,
\begin{equation*}
f(q)=\left\vert dR_{B_{1}^{-1}}dL_{A_{1}}\left( \left. \frac{d}{dt}%
\right\vert _{0}\gamma (t)\ast q\right) \right\vert ^{2}=\left\vert \left.
\frac{d}{dt}\right\vert _{0}\gamma (t)\ast q\right\vert ^{2}\text{,}
\end{equation*}
and so
\begin{equation*}
\left\Vert dF_{k_{1},k_{2}}X\right\Vert ^{2}=\int_{M}\left\vert \left. \frac{%
d}{dt}\right\vert _{0}\gamma (t)\ast q\right\vert ^{2}d\mu(q) =\left\Vert
X\right\Vert ^{2}\text{.}
\end{equation*}%
Then $F_{k_{1},k_{2}}$ is an isometry.}

\medskip
\textit{
b) We start with the least simple case $\mathbb{F}=\mathbb{C}$, when $K$ is
strictly contained in $M\times M$, that is, $K=S\left( U_{n}\times
U_{n}\right) $ and $M=U_{n}$. The cases $M=SO_{n}$ (with $n>2$) and $Sp_{n}$%
, which are semisimple Lie groups, are dealt with in the same way, but in a
more direct manner. We have to show that the kinetic energy metric
restricted to $K=S\left( U_{n}\times U_{n}\right) $ is the restriction of a
product metric on $M\times M=U_{n}\times U_{n}$ such that $\rho$ is an
isometry.}

\textit{
By the first part of the proof, the kinetic energy metric restricted to $K$
is bi-invariant, in particular, the inner product on $\mathfrak{k}%
=T_{\left(  I_{n},I_{n}\right)}K$ is Ad$\left( K\right) $-invariant. We write
\begin{equation*}
\mathfrak{k}=s\left( u_{n}\oplus u_{n}\right) =su_{n}\oplus su_{n}\oplus
\mathbb{R}iI_{2n}^{-}\text{,}
\end{equation*}%
where $I_{2n}^{-}=\operatorname{diag}\left( I_{n},-I_{n}\right) $. Since $su_{n}$ is
simple and the last direct summand has dimension 1, the summands are
orthogonal.}

\textit{
Now we check that $d\rho _{\left( I_{n},I_{n}\right) }$ restricted to $%
su_{n}\oplus su_{n}$ is a linear isometry, or equivalently, that $\left\Vert
\left( X,0\right) \right\Vert =\left\Vert \left( 0,X\right) \right\Vert $
for all $\left( X,0\right) \in \mathfrak{k}$ (using again that $\mathfrak{k}$
carries the restriction of a product metric). We compute
\begin{equation*}
\left\vert \left. \frac{d}{dt}\right\vert _{0}e^{t\left( X,0\right) }\ast
q\right\vert ^{2}=\left\vert \left. \frac{d}{dt}\right\vert
_{0}e^{tX}q\right\vert ^{2}=\left\vert Xq\right\vert ^{2}=\left\vert
X\right\vert ^{2}\text{,}
\end{equation*}%
for all $q\in M$, since the inner product on $\mathbb{F}^{m\times m}$ is
bi-invariant by $M.$ Integrating over $M$ yields $\left\Vert \left(
X,0\right) \right\Vert ^{2}=\operatorname{vol}\left( M\right) \left\vert
X\right\vert ^{2}$. One obtains similarly that $\left\Vert \left( 0,X\right)
\right\Vert ^{2}=\operatorname{vol}\left( M\right) \left\vert X\right\vert ^{2}$.}

\textit{
Finally, we verify that the metric on $K$ is the restriction to $K$ of a
product metric on $M\times M$. We observe that $u_{n}\oplus u_{n}=\mathfrak{k%
}\oplus \mathbb{R}iI_{2n}$. Let $\epsilon _{1}=\operatorname{diag}\left(
iI_{n},0_{n}\right) $ and $\epsilon _{2}=\operatorname{diag}\left(
0_{n},iI_{n}\right) $, which form a basis of $\mathbb{R}iI_{2n}^{-}\oplus
\mathbb{R}iI_{2n}$. We have that
\begin{equation*}
u_{n}\oplus u_{n}=\left( su_{n}\oplus \mathbb{R}\epsilon _{1}\right) \oplus
\left( su_{n}\oplus \mathbb{R}\epsilon _{2}\right) \equiv \left(
su_{n}\oplus su_{n}\right) \oplus \mathbb{R}\epsilon _{1}\oplus \mathbb{R}%
\epsilon _{2}\text{.}
\end{equation*}%
Suppose that $\left\Vert iI_{2n}^{-}\right\Vert =\ell $. We consider an
inner product such that the three summands in the last expression are
orthogonal, $su_{n}\oplus su_{n}$ has the metric induced from the kinetic
energy and $\left\Vert \epsilon _{1}\right\Vert =\left\Vert \epsilon
_{2}\right\Vert =\ell /\sqrt{2}$. Since $iI_{2n}^{-}=\epsilon _{1}-\epsilon
_{2}$, this inner product on $u_{n}\oplus u_{n}$ satisfies the required
conditions.}\qed

\section{Some totally geodesic subgroups of $G$}

Before proving Theorem \ref{TotGeod} we comment on the canonical inclusions.
First, notice that
\begin{equation}
U\left( n,n\right) =\{A\in Sp(n,n)\mid iA=Ai\}\text{, \ } O\left( n,n\right)
=\{A\in Sp(n,n)\mid iA=Ai\text{ y }jA=Aj\}\text{,}  \label{OnnCSpnn}
\end{equation}
since for each entry $q=a+bi+cj+dk$ of a matrix in $Sp\left( n,n\right) $, $%
iq=qi$ holds if and only if $c=d=0$, and also $jq=qj$ if and only if $b=d=0$.

Now we realize the inclusions of $U(n,n)$ and $Sp\left( n,n\right) $ as
subgroups of $O(2n,2n)$ and $O(4n,4n)$, respectively. We follow Chapter 2 of
\cite{Har}. We recall first that under a suitable identification of $\mathbb{%
C}^{2n}$ with $\mathbb{R}^{4n}$, it is possible to describe%
\begin{equation}
U(n,n)=\{A\in O_{o}(2n,2n)\mid JA=AJ\},  \label{Unn}
\end{equation}%
where $J:\mathbb{R}^{4n}\rightarrow \mathbb{R}^{4n}$ is the matrix with $2n$
diagonal blocks of the form
$J_{0}=\left(
\begin{array}{cc}
0 & -1 \\
1 & 0%
\end{array}%
\right) $.

Consider now the right vector space ${\mathbb{H}}^{2n}$ over $\mathbb{H}$.
Now, any $q\in {\mathbb{H}}$ has the form $q=a+bi+cj+dk$, with $a,b,c,d\in {%
\mathbb{R}}$ and so ${\mathbb{H}}^{2n}$ is naturally isomorphic to ${\mathbb{%
R}}^{8n}$. Via this identification,
\begin{equation}
Sp\left( n,n\right) =\{A\in O_{o}(4n,4n)\mid R_{i}A=AR_{i}\text{ and }%
R_{j}A=AR_{j}\}\text{,}  \label{Spnn}
\end{equation}%
where $R_{p}$ denotes right multiplication by $p\in {\mathbb{H}}$. The ($%
8n\times 8n$)-matrices of $R_{i}$ and $R_{j}$ with respect to the canonical
basis consist, respectively, of $2n$ diagonal blocks of the forms
\begin{equation*}
\left(
\begin{array}{cc}
J_{0} & 0 \\
0 & -J_{0}%
\end{array}%
\right) \text{\ \ \ \ \ \ and\ \ \ \ \ \ }\left(
\begin{array}{cc}
0 & -I_{2} \\
I_{2} & 0%
\end{array}%
\right) \text{.}
\end{equation*}

The following lemma will be useful to prove that $O_{o}\left( n,n\right) $
is totally geodesic in $SU\left( n,n\right) $. Observe that the arguments in
the proof are not valid for $Sp(n,n)$ (see \cite{Zh}).

\begin{lemn}
\label{Conj} If $G=SU\left( n,n\right) $ is endowed with the kinetic energy
metric, then conjugation $C:G\rightarrow G$, $C\left( q\right) =\overline{q}%
, $ is an isometry.
\end{lemn}

\begin{proofn}
Let $X\in T_{p}G$ and let $\gamma $ be a curve such that $\gamma (0)=p$ and $%
\gamma ^{\prime }(0)=X$. We have%
\begin{eqnarray*}
\Vert dC_{p}(X)\Vert ^{2} &=&\int_{U_{n}}\left\vert \left. \frac{d}{dt}%
\right\vert _{0}\left( \overline{\gamma (t)}\ast q\right) \right\vert
^{2}d\mu =\int_{U_{n}}\left\vert \overline{\left. \frac{d}{dt}\right\vert
_{0}(\gamma (t)\ast \overline{q})}\right\vert ^{2}d\mu \\
&=&\int_{U_{n}}\left\vert \left. \frac{d}{dt}\right\vert _{0}(\gamma (t)\ast
\overline{q})\right\vert ^{2}d\mu =\Vert X\Vert ^{2}\text{.}
\end{eqnarray*}%
The last equality holds by the change of variables theorem, since
conjugation in $U_n$ is an isometry of $U_n$. \qed
\end{proofn}

\medskip
\noindent
\textbf{Proof of Theorem \protect\ref{TotGeod}}.\textit{
We use repeatedly that each connected component of the set of fixed points
of an isometry is a totally geodesic submanifold (see for instance
Proposition 10.3.6 in \cite{Olmos}). This and the previous lemma yield that
the the first inclusion is totally geodesic, since $O_{o}(n,n)$ is a
connected component of the set of fixed points of the conjugation.}

\textit{
For $k\in K$, let $F_{k}:G\rightarrow G$ be the conjugation by $k$, that is, 
$F_{k}\left( g\right) =kgk^{-1}$. Notice that $F_{k}=F_{k,k}$ with the
notation of (\ref{Fk1k2}). By Proposition \ref{CompactoMaximal} (a), it
suffices to show in each case the existence of $k\in K$ such that the
corresponding subgroup is the identity component of the set of fixed points
of the isometry $F_{k}$. Actually, in the cases (c) and (e) we need to
consider the intersection of such sets for different elements of $K$. It is
not difficult to see that the following assertion is true: Let $N_{1}$ and $%
N_{2}$ be two totally geodesic submanifolds. If a connected component of $%
N_{1}\cap N_{2}$ is a submanifold, then it is totally geodesic.}

\textit{
Now, $U\left( n,n\right) $ is a totally geodesic subgroup of $Sp\left(
n,n\right) $ since it is connected and by (\ref{OnnCSpnn}) is the set of
fixed points of $F_{iI_{2n}}$, which is an isometry of $Sp\left( n,n\right) $%
, as $iI_{n}\in Sp_{n}$. Analogously, $O_{o}\left( n,n\right) $ is a totally
geodesic subgroup of $Sp\left( n,n\right) $, given that (\ref{OnnCSpnn}) is
the identity component of the intersection of the sets of fixed points of $%
F_{iI_{2n}}$ and $F_{jI_{2n}}$, which are isometries of $Sp\left( n,n\right) 
$.}

\textit{
By (\ref{Unn}), the group $U\left( n,n\right) $ is the set of fixed points
of the map $F_{J}:O_{o}(2n,2n)\rightarrow O_{o}(2n,2n)$ given by $%
F_{J}(A)=JAJ^{-1}=-JAJ$. Now, by the description of $J$ above, $J\in
SO_{2n}\times SO_{2n}$. Then, $F_{J}$ is an isometry of $O_{o}(2n,2n)$
endowed with the kinetic energy metric. Since $U(n,n)$ is connected, it is a
totally geodesic submanifold of $O_{o}(2n,2n)$.}

\textit{
By (\ref{Spnn}), the group $Sp\left( n,n\right) $ is the intersection of the
sets of fixed points of the maps $F_{R_{i}},F_{R_{j}}:O_{o}(4n,4n)%
\rightarrow O_{o}(4n,4n)$ given by $F_{R_{i}}\left( A\right)
=R_{i}AR_{i}^{-1}$, and similarly for $F_{R_{j}}$. From the form of the
matrices of $R_{i}$ and $R_{j}$ (see above) one deduces that they belong to $%
SO_{4n}\times SO_{4n}$. Thus, $F_{R_{i}}$ and $F_{R_{j}}$ are isometries of $%
O_{o}(4n,4n)$ endowed with the kinetic energy metric. Thus, $Sp\left(
n,n\right) $, which is connected, is totally geodesic in $O_{o}(4n,4n)$.}
\qed

\section{The incompleteness of the kinetic energy metric}

In this section we prove Theorem \ref{geodfinita copy(1)}. Let $\delta =%
\operatorname{diag}(I_{n},-I_{n})$. Then $b_{n,n}^{\mathbb{F}}$ as in (\ref{form})
is given by $b_{n,n}^{\mathbb{F}}(x,y)=\overline{x}^{T}\delta y$, for $%
x,y\in \mathbb{F}^{2n}$. Thus, $U\in U(n,n,\mathbb{F})$ if and only if $%
\overline{U}^{T}\delta U=\delta $. From this expression one gets that $u\in
\mathfrak{u}(n,n,\mathbb{F})$, the Lie algebra of $U(n,n,\mathbb{F})$, if
and only if $\overline{u}^{T}\delta +\delta u=0$. Hence,
\begin{equation}
\mathfrak{g}=\left\{ \left(
\begin{array}{cc}
a & c \\
\overline{c}^{T} & b%
\end{array}%
\right) \mid c\in \mathbb{F}^{n\times n}\text{, }\overline{a}^{T}=-a\text{, }%
\overline{b}^{T}=-b\right\}  \label{LieG}
\end{equation}%
is the Lie algebra of $G$ for $\mathbb{F}=\mathbb{R}$ and $\mathbb{H}$,
whereas the elements of $su\left( n,n\right) $ have the same form, with the
extra condition $\operatorname{tr}(a+b)=0$.

\begin{lemn}
\label{lemapuntosfijos}For $i=2,\dots ,n$, let $T_{i}$ be the $\mathbb{F}$%
-linear transformation on $\mathbb{F}^{n}$ such that $T_{i}\left(
e_{1}\right) =e_{i}$, $T_{i}\left( e_{i}\right) =-e_{1}$ and $T_{i}\left(
e_{s}\right) =e_{s}$ for $1<s\leq n$, $s\neq i$. Then an $\mathbb{F}$-linear
transformation $T$ on $\mathbb{F}^{n}$ commutes with $T_{i}$ for all $%
i=2,\dots ,n$ if and only if $T=\lambda I_{n}$, for some $\lambda \in
\mathbb{F}$.
\end{lemn}

\begin{proofn}
The matrix of $T_{i}$ with respect to the canonical basis is $A^{i}$, where
\begin{equation*}
A_{1,i}^{i}=1\text{,\ \ \ }A_{i,1}^{i}=-1\text{, \ }A_{j,j}^{i}=1\text{ \ \
for }1<j\neq i
\end{equation*}%
and $A_{s,t}^{i}=0$ otherwise. {First observe that for any matrix $X$, $%
A^{i}X$ is the matrix whose rows coincide with those of }$X$, except that
its {first row is the $i$-th row of $X$ and its $i$-th row is minus the
first row of $X$. On the other hand, $XA^{i}$ is obtained by a similar
procedure, but interchanging columns of $X$ instead of rows. So we obtain
that }
\begin{equation}
X_{ii}=\left( A^{i}X\right) _{1,i}=\left( XA^{i}\right) _{1,i}=X_{11}\text{.}
\label{cond1}
\end{equation}%
for all $1<i\leq n$. Also, for $t\neq i$,
\begin{equation}
X_{ti}=\left( A^{i}X\right) _{ti}=\left( XA^{i}\right) _{ti}=X_{t1}=\left(
A^{i}X\right) _{t1}=\left( XA^{i}\right) _{t1}=-X_{ti}  \label{cond2}
\end{equation}

Combining conditions (\ref{cond1}) and (\ref{cond2}), we get $X_{ii}=X_{11}$
for every $i$ and $X_{ti}=0$ for $t\neq i$. Therefore $X=\lambda I_{n}$, for
some $\lambda \in \mathbb{F}$.
\qed

\smallskip

\medskip
\noindent
\textbf{\emph{ Proof of Theorem \protect\ref{geodfinita copy(1)}}.}
The proof will consist of two parts.

\smallskip

\textbf{Part 1. }We prove that the image of $\gamma $ is a totally geodesic
submanifold of $G$ and so $\gamma $ is the reparametrization of a geodesic
in $G$, which will be inextendible since $\lim_{t\rightarrow \pm \infty
}\gamma \left( t\right) $ does not exist.

For $i=2,\dots ,n$, let $k_{i}=\operatorname{diag}(A^{i},A^{i})\in K$, where $A^{i}$
is the matrix introduced in Lemma \ref{lemapuntosfijos} (notice that $A^i\in
SO_n$), and let $F_{i}:G\rightarrow G,$ $F_{i}\left( g\right)
=k_{i}gk_{i}^{-1}$. Let $L$ be the identity component of the Lie group 
\begin{equation*}
	\left( \cap _{i=2}^{n}\,\operatorname{Fix}_{0}\left( F_{i}\right) \right) \cap
	O_{o}\left( n,n\right) \text{.}
\end{equation*}

We will prove that $L$ coincides with the image of $\gamma $. Thus, this
will be totally geodesic in $G$, since Fix$_{0}\left( F_{i}\right) $ and $%
O_{o}\left( n,n\right) $ are totally geodesic in $G$ by Proposition \ref%
{CompactoMaximal} (a) and Theorem \ref{TotGeod} (a) and (c), respectively
(as we mentioned in the proof of the latter, if a connected component of the
intersection of a finite number of totally geodesic submanifolds is a
submanifold, then it is totally geodesic).

The Lie algebra of $L$ is $\mathfrak{f}\cap o\left( n,n\right) $, where 
\begin{equation*}
	\mathfrak{f}=\cap _{i=2}^{n}\,\operatorname{Fix}\left( \operatorname{Ad}\left( k_{i}\right)
	\right) \text{,}
\end{equation*}%
since the Lie algebra of the identity component of the intersection of a
finite number of Lie subgroups coincides with the intersection of its Lie
algebras.

For each $i=2,\dots,n$ we apply $\operatorname{Ad}(k_{i})$ to a generic element of $%
\mathfrak{g}$ as in the beginning of the section. One easily verifies that 
\begin{equation*}
	\operatorname{Ad}(k_{i})\left( 
	\begin{array}{cc}
		a & \overline{c}^{T} \\ 
		c & b%
	\end{array}
	\right) =\left( 
	\begin{array}{cc}
		A_{i}aA_{i}^{-1} & A_{i}\overline{c}^{T}A_{i}^{-1} \\ 
		A_{i}cA_{i}^{-1} & A_{i}bA_{i}^{-1}%
	\end{array}
	\right) \text{.}
\end{equation*}

Since $\operatorname{Fix}\left( \operatorname{Ad}\left( k_{i}\right) \right) =\left\{ X\in%
\mathfrak{g}\mid X\text{ conmuta con }k_{i}\right\}$, by the conditions on $%
a $ and $b$ in (\ref{LieG}), using Lemma \ref{lemapuntosfijos} we have that 
\begin{equation*}
	\mathfrak{f}=\left\{ \left( 
	\begin{array}{cc}
		\alpha I_{n} & \overline{\lambda}I_{n} \\ 
		\lambda I_{n} & \beta I_{n}%
	\end{array}
	\right) \mid\alpha,\lambda\in\mathbb{F}\text{, }\overline{\alpha}=-\alpha%
	\text{, }\overline{\beta}=-\beta\right\}
\end{equation*}
(requiring additionally that $\alpha+\beta=0$ in the complex case). Now,
since $o\left( n,n\right) $ consists of matrices whose entries are real
numbers, we have that 
\begin{equation*}
	\mathfrak{f}\cap o\left( n,n\right) =\mathbb{R}\left( 
	\begin{array}{cc}
		0 & I_{n} \\ 
		I_{n} & 0%
	\end{array}
	\right) \text{.}
\end{equation*}
Therefore, the Lie algebras of $L$ and $\gamma$ coincide. Consequently, $%
L=\gamma\left( \mathbb{R}\right) $, given that both are connected.

\medskip

\textbf{Part 2. }We prove that $\gamma $ has finite length. It will be
convenient to consider the reparametrization $\sigma :\left( -1,1\right)
\rightarrow G$ of $\gamma $ given by 
\begin{equation*}
	\sigma (s)=\gamma (\tanh ^{-1}(s))=\frac{1}{\sqrt{1-s^{2}}}\left( 
	\begin{array}{cc}
		I_{n} & sI_{n} \\ 
		sI_{n} & I_{n}%
	\end{array}%
	\right) \text{.}
\end{equation*}%
Then, by (\ref{dinmet}), 
\begin{equation}
	\left\Vert \sigma ^{\prime }\left( t\right) \right\Vert
	^{2}=\int_{M}\left\vert \widetilde{\sigma ^{\prime }(t)}(A)\right\vert
	^{2}d\mu \left( A\right) =\int_{M}\left\vert \left. \frac{d}{ds}\right\vert
	_{t}\sigma (s)\ast A\right\vert ^{2}d\mu \left( A\right) .
	\label{normaSigmaP}
\end{equation}

Let $T=\left\{ D(\theta _{1},\cdots ,\theta _{m})\mid \theta _{j}\in \mathbb{%
	R}\right\} $ be a maximal torus of $M$, where $D(\theta _{1},\cdots ,\theta
_{m})$ equals

- $\operatorname{diag}(R_{\theta _{1}},\cdots ,R_{\theta _{m}})$ with $R_{\theta
}=\left( 
\begin{array}{cc}
	\cos \theta & -\sin \theta \\ 
	\sin \theta & \cos \theta%
\end{array}%
\right) $, if $M=SO_{2m}$.

- $\operatorname{diag}(R_{\theta _{1}},\cdots ,R_{\theta _{m}},1)$, if $M=SO_{2m+1}$.

\smallskip

- $\operatorname{diag}(e^{i\theta _{1}},\cdots ,e^{i\theta _{m}})$, if $M=U_{m}$ or $%
M=Sp_{m}$. 
\smallskip

We recall the general form of the Weyl integration formula on $M$. Provided
that $f$ is a class function, that is, $f\left( BAB^{-1}\right) =f\left(
A\right) $ for all $A,B\in M$, we have 
\begin{equation*}
	\int_{M}f\left( A\right) d\mu \left( A\right) =\int_{T}f\left( D(\theta
	_{1},\cdots ,\theta _{m})\right) F\left( \theta _{1},\cdots ,\theta
	_{m}\right) d\theta _{1}\cdots d\theta _{m}\text{,}
\end{equation*}%
where $F$ is some bounded function. The precise formula for $U_{n}$ can be
found in Theorem 3.1 in \cite{Me}, and for the remaining groups, in Theorems
IX.9.3, IX.9.4, IX.9.5 in \cite{Si}.

It is not difficult to see that for every{\ $A,B\in M$, $\sigma (s)\ast
	(BAB^{-1})=B(\sigma (s)\ast A)B^{-1}$ and this implies that the integrand in
	the right hand side of (\ref{normaSigmaP}) is a class function of }$A$.
Then, 
\begin{equation}
	\int_{M}\left\vert \left. \frac{d}{ds}\right\vert _{t}\sigma (s)\ast
	A\right\vert ^{2}d\mu \left( A\right) \leq C\int_{T}\left\vert \left. \frac{d%
	}{ds}\right\vert _{t}\sigma (s)\ast D(\theta _{1},\cdots ,\theta
	_{m})\right\vert ^{2}d\theta _{1}\cdots d\theta _{m}  \label{Weyl}
\end{equation}
for some constant $C>0$.

Consider first the case $M=SO_{n}$ with $n=2m$. Since 
\begin{equation*}
	\sigma (s)\ast D(\theta _{1},\cdots ,\theta _{m})=\operatorname{diag}(\sigma
	_{2}(s)\ast R_{\theta _{1}},\cdots ,\sigma _{2}(s)\ast R_{\theta _{m}}),
\end{equation*}%
where $\sigma _{2}(s)=\frac{1}{\sqrt{1-s^{2}}}\left( 
\begin{array}{cc}
	I_{2} & sI_{2} \\ 
	sI_{2} & I_{2}%
\end{array}%
\right) \in O_{o}\left( 2,2\right) $, then 
\begin{equation*}
	\left\vert \left. \frac{d}{ds}\right\vert _{t}\sigma (s)\ast D(\theta
	_{1},\cdots ,\theta _{m})\right\vert ^{2}=\sum_{j=1}^{m}\left\vert \left. 
	\frac{d}{ds}\right\vert _{t}\sigma _{2}(s)\ast R_{\theta _{j}}\right\vert
	^{2}\text{.}
\end{equation*}%
We compute 
\begin{eqnarray}
	\sigma _{2}\left( s\right) \ast R_{\theta } &=&\left( R_{\theta
	}+sI_{2}\right) \left( sR_{\theta }+I_{2}\right) ^{-1}  \notag \\
	&=&\left( 
	\begin{array}{cc}
		s+\cos \theta & -\sin \theta \\ 
		\sin \theta & s+\cos \theta%
	\end{array}%
	\right) \left( 
	\begin{array}{cc}
		1+s\cos \theta & -s\sin \theta \\ 
		s\sin \theta & 1+s\cos \theta%
	\end{array}%
	\right) ^{-1}  \notag \\
	&=&\frac{1}{1+2s\cos \theta +s^{2}}\left( 
	\begin{array}{cc}
		2s+\left( 1+s^{2}\right) \cos \theta & (s^{2}-1)\sin \theta _{j} \\ 
		(1-s^{2})\sin \theta & 2s+\left( 1+s^{2}\right) \cos \theta%
	\end{array}%
	\right) \text{.}  \label{sigma2}
\end{eqnarray}

Now we identify $SO_{2}$ with $S^{1}$ in the usual way, that is, $R_{\alpha
} $ is identified with $e^{i\alpha }$. Then (\ref{sigma2}) reads%
\begin{equation*}
	\sigma _{2}\left( s\right) \ast e^{i\theta }=\frac{e^{i\theta }+s}{%
		1+se^{i\theta }}\text{,}
\end{equation*}%
since the real and imaginary parts of this complex number are equal to the
first and second entries of the first column of the matrix (\ref{sigma2}),
respectively. Replacing in (\ref{Weyl}) we obtain 
\begin{equation*}
	\left\Vert \sigma ^{\prime }\left( t\right) \right\Vert
	^{2}=\int_{M}\left\vert \left. \frac{d}{ds}\right\vert _{t}\sigma (s)\ast
	A\right\vert ^{2}d\mu \left( A\right) \leq
	C\sum_{j=1}^{m}\int_{S^{1}}\left\vert \left. \frac{d}{ds}\right\vert _{t}%
	\frac{e^{i\theta _{j}}+s}{1+se^{i\theta _{j}}}\right\vert ^{2}d\theta _{j}%
	\text{.}
\end{equation*}

Setting $z=e^{i\theta _{j}}$, the integrals in the right hand side can be
written as integrals over the circle of a meromorphic function on the
two-punctured disc. They can be computed by residues and shown to be equal
to $4\pi /\left( 1-t^{2}\right) $. This is done in \cite{ES}, where the
force free M\"{o}bius motions of the circle are studied (the $j$-th
integrand equals $\left\vert \widetilde{Y}\left( e^{i\theta _{j}}\right)
\right\vert ^{2}$, with $\widetilde{Y}$ as in the proof of Theorem 3 of that
article). Therefore,%
\begin{equation*}
	\text{length}\left( \left. \sigma \right\vert _{\left[ 0,1\right) }\right)
	=\lim_{T\rightarrow 1^{-}}\int_{0}^{T}\left\Vert \sigma ^{\prime }\left(
	t\right) \right\Vert ~dt\leq C^{\prime }\int_{0}^{1}\frac{2}{\sqrt{1-t^{2}}}%
	\,dt=C^{\prime }\pi
\end{equation*}
(with $C^{\prime }=\sqrt{Cm\pi }$) and so the length of $\left. \sigma
\right\vert _{\left[ 0,1\right) }$ is finite, as desired.

The case $M=SO_{2m+1}$ is analogous. If $M=U_{n}$ or $Sp_{n}$, we get that 
\begin{eqnarray*}
	\sigma (s)\ast D(\theta _{1},\cdots ,\theta _{n}) &=&(D(\theta _{1},\cdots
	,\theta _{n})+sI_{n})(sD(\theta _{1},\cdots ,\theta _{n})+I_{n})^{-1} \\
	&=&\operatorname{diag}\left( \frac{e^{i\theta _{1}}+s}{1+se^{i\theta _{1}}},\cdots ,%
	\frac{e^{i\theta _{n}}+s}{1+se^{i\theta _{n}}}\right)
\end{eqnarray*}%
and the proof follows from similar arguments.
\qed

\medskip
\noindent
\textbf{\emph{Proof of Proposition \protect\ref{acum}}.}
For $q\in M$ we compute%
\begin{equation*}
\gamma \left( t\right) \ast q =\left( \cosh t\ q+\sinh t\ I_{n}\right)
\left( \sinh t\ q+\cosh t\ I_{n}\right) ^{-1}= \left( q+\tanh t\
I_{n}\right) \left( \tanh t\ q+I_{n}\right) ^{-1}\text{. }
\end{equation*}

For $\varepsilon =\pm 1$, if $-\varepsilon $ is not an eigenvalue of $q$, we
have that%
\begin{equation*}
\lim_{t\rightarrow \varepsilon \infty }\gamma \left( t\right) \ast q=\left(
q+\varepsilon I_{n}\right) \left( \varepsilon q+I_{n}\right)
^{-1}=\varepsilon I_{n}\text{.}
\end{equation*}

Suppose that $M=SO_{n}$ with $n$ even. Given $q\in M$ with eigenvalue $%
-\varepsilon ,$ we verify that any neighborhood of $q$ contains elements
with no eigenvalue $-\varepsilon .$ Let $R_{t}$ denote the rotation in $%
\mathbb{R}^{2}$ through the angle $t$. We have
\begin{equation*}
q=O\operatorname{diag}\left( R_{\alpha _{1}},\dots ,R_{\alpha _{k}},-\varepsilon
I_{\ell }\right) O^{-1}
\end{equation*}%
for some $O\in SO_{n}$, where $\ell $ is even and $R_{\alpha _{i}}$ has no
eigenvalue $-\varepsilon $. Now, the curve%
\begin{equation*}
t\mapsto q_{t}:=O\operatorname{diag}\left( R_{\alpha _{1}},\dots ,R_{\alpha
_{k}},-\varepsilon R_{t},\dots ,-\varepsilon R_{t}\right) O^{-1}
\end{equation*}%
has initial point $q$ and $q_{t}$ has no eigenvalue $-\varepsilon $ for $t>0$
near $0$. Similar arguments hold for $SO_{n}$ with $n$ odd if $\varepsilon
=1 $.

The cases $M=U_{n}$ or $Sp_{n}$ can be dealt with analogously, with $%
e^{i\alpha _{j}}$ and and $e^{it}$ instead of $R_{\alpha _{j}}$ and $R_{t}$,
respectively.
\qed
\end{proofn}

\section{Force free rigid motions}


We recall that the maximal compact subgroup $K$ of $G$ is a product $%
SO_{n}\times SO_{n}$ or $Sp_{n}\times Sp_{n}$ or is contained as the
hypersurface $S\left( U_{n}\times U_{n}\right) $ in the product $U_{n}\times
U_{n}$. Theorem \ref{XceroceroY} asserts that geodesics in $K$ contained in
one of the factors are also geodesics in the ambient Lie group $G$ endowed
with the kinetic energy metric. Next we prove the theorem.

\medskip
\noindent
\textbf{Proof of Theorem \protect\ref{XceroceroY}}.
\emph{
Let $\widehat{Z}$ be the right invariant vector field of $G$ such that $%
\widehat{Z}_{I_{2n}}=Z$. Then its flow is $t\mapsto\varphi _{t}=L_{\exp
(tZ)} $. Now, for all $t$, $\varphi _{t}$ is an isometry of $G$ by
Proposition \ref{CompactoMaximal} (a). Therefore, $\widehat{Z} $ is a
Killing vector field of $G$ and it is standard to check that $\gamma $ is an
integral curve of $\widehat{Z}$ through the identity. So, in order to see
that $\gamma $ is a geodesic of $G$, we need to show that $(\nabla _{%
\widehat{Z}}\widehat{Z})_{\gamma (t)}=0$, where $\nabla $ is the Levi-Civita
connection of $G$ endowed with the kinetic energy metric. Since $K$ acts on $%
G$ by isometries and $dR_{\exp \left( tZ\right) }(\widehat{Z})=\widehat{Z}\circ R_{\exp \left( tZ\right)}$,
it suffices to see that $(\nabla _{\widehat{Z}}\widehat{Z})_{I_{2n}}=0$, or
equivalently, since $\widehat{Z}$ is Killing, that
\begin{equation}
-2\langle \nabla _{Z}\widehat{Z},Y\rangle =2\langle \nabla _{Y}\widehat{Z}%
,Z\rangle =Y(\Vert \widehat{Z}\Vert ^{2})=\left. \frac{d}{dt}\right\vert
_{0}\Vert \widehat{Z}(\exp (tY))\Vert ^{2}  \label{conexion1}
\end{equation}%
vanishes for every $Y\in \mathfrak{g}$ (the second equality holds by the
Koszul formula). We have that
\begin{equation*}
\widehat{Z}(\exp (tY))=\left. \frac{d}{ds}\right\vert _{0}e^{sZ}e^{tY}
\end{equation*}
and so, calling $q_{t}=e^{tY}\ast q\in M$, by (\ref{dinmet}),
\begin{equation}
\Vert \widehat{Z}(\exp (tY))\Vert ^{2}=\int_{M}\left\vert \left. \frac{d}{ds}%
\right\vert _{0}\left( e^{sZ}e^{tY}\right) \ast q\right\vert ^{2}d\mu \left(
q\right) =\int_{M}\left\vert \left. \frac{d}{ds}\right\vert _{0}e^{sZ}\ast
q_{t}\right\vert ^{2}d\mu \left( q\right)\text{.}  \label{integral}
\end{equation}}

\emph{
Now suppose that $Z=\operatorname{diag}\left( X,Y\right) $ with $X,Y\in \mathfrak{m}$
and $X=0$ or $Y=0$. Hence, by (\ref{mainAction}),
\begin{equation*}
e^{sZ}\ast q_{t}=\operatorname{diag}\left( e^{sX},e^{sY}\right) \ast
q_{t}=e^{sX}q_{t}e^{-sY}
\end{equation*}%
and so the integrand in (\ref{integral}) is $\left\vert
Xq_{t}-q_{t}Y\right\vert ^{2}$, which equals $\left\vert Xq_{t}\right\vert
^{2}=\left\vert X\right\vert ^{2}$ or $\left\vert q_{t}Y\right\vert
^{2}=\left\vert Y\right\vert ^{2}$ if $Y=0$ or $X=0$, respectively (notice
that the inner product on $\mathbb{F}^{n\times n}$ is bi-invariant). In both
cases we have independence on $t$. Therefore, the expression (\ref{conexion1}%
) vanishes, as desired.}
\qed

\section{Equivalences with conformal and projective motions of the sphere in
low dimensions \label{last}}

As we mentioned in the introduction, in this section we explore, on the one
hand, whether the projective and conformal actions on spheres (\ref%
{projective}) are equivalent to the M\"{o}bius actions on unitary groups in
low dimensions, and on the other hand, the implications of this in the
dynamics of M\"obius motions of the compact classical Lie groups.

The simplest case is the M\"{o}bius action $SU\left( 1,1\right) \times
U_{1}\rightarrow U_{1}$. The equivalences up to coverings with the usual
actions of $SL_{2}\left( \mathbb{R}\right) $ and $O_{o}\left( 1,2\right) $
on the circle (M\"{o}bius or projective and conformal, respectively) is well
known. In \cite{ES} the first two authors gave a complete description of the
geometry of the kinetic energy metric in this simple case. See also \cite%
{ESarXiv}, where the relationship with the article \cite{BetterThanNice} is
addressed, providing a dynamical presentation of the metric on the disc
where everything is better than nice.\ The kinetic energy metric projects to a non-complete negatively curved disc
whose geodesics are hypocycloids.

The next simplest case is the M\"{o}bius action $O_{o}\left(
2,2\right) \times SO_{2}\rightarrow SO_{2}$, which is the only one of those
actions which is not even almost effective. For the sake of completeness we
comment on its relationship with projective transformations of the circle.

We recall from \cite{Har} the split quaternions $\mathbb{M}$, which is the
skew algebra $M\left( 2,\mathbb{R}\right) =\mathbb{R}^{2\times 2}$ with the
usual matrix multiplication and the inner product whose square norm is given
by $\left\langle X,X\right\rangle ={\mathrm{det}}X$. In particular, $%
SL\left( 2,\mathbb{R}\right) $ is the unit sphere in $\mathbb{M}$. The
matrices
\begin{equation*}
\mathbf{1}=\left(
\begin{array}{cc}
1 & 0 \\
0 & 1%
\end{array}%
\right) \text{,\ \ \ \ }\mathbf{i}=\left(
\begin{array}{cc}
0 & 1 \\
-1 & 0%
\end{array}%
\right) \text{,\ \ \ \ }\mathbf{j}=\left(
\begin{array}{cc}
0 & 1 \\
1 & 0%
\end{array}%
\right) \text{,\ \ \ \ }\mathbf{k}=\left(
\begin{array}{cc}
1 & 0 \\
0 & -1%
\end{array}%
\right)
\end{equation*}%
form a basis of $\mathbb{M}$ and the following relations hold:%
\begin{equation*}
-\mathbf{i}^{2}=\mathbf{j}^{2}=\mathbf{k}^{2}=\mathbf{1}\text{,\ \ \ \ }%
\mathbf{ij}=-\mathbf{ji}=\mathbf{k}\text{,\ \ \ \ }\mathbf{jk}=-\mathbf{kj}=-%
\mathbf{i}\text{ \ \ \ and\ \ \ \ }\mathbf{ki}=-\mathbf{ik}=\mathbf{j}\text{.%
}
\end{equation*}%
The square norm is given by
\begin{equation*}
\left\langle a+b\mathbf{i}+c\mathbf{j}+d\mathbf{k},a+b\mathbf{i}+c\mathbf{j}%
+d\mathbf{k}\right\rangle =a^{2}+b^{2}-c^{2}-d^{2}
\end{equation*}%
and hence the inner product has signature $\left( 2,2\right) $. Observe that
$p\overline{p}=\left\langle p,p\right\rangle $ for all $p\in \mathbb{M}$ and
$p^{-1}=\overline{p}/\left\langle p,p\right\rangle $ if $\left\langle
p,p\right\rangle \neq 0$. Let $L_{p}$ and $R_{p}$ denote left and right
multiplication by $p$, respectively, which are linear isometries of $\mathbb{%
M}$ if $\left\langle p,p\right\rangle =1$. One has that%
\begin{equation*}
SL\left( 2,\mathbb{R}\right) \times SL\left( 2,\mathbb{R}\right)
\longrightarrow O_{o}\left( \mathbb{M}\right) \text{,}\ \ \ \ \ \ \ \left(
p,q\right) \mapsto L_{p}\circ R_{\overline{q}}
\end{equation*}%
is a double covering morphism with kernel $\left\{ \left( \mathbf{1},\mathbf{%
1}\right) ,\left( -\mathbf{1},-\mathbf{1}\right) \right\} $. In particular,%
\begin{equation*}
O_{o}\left( \mathbb{M}\right) =\left\{ L_{p}\circ R_{\overline{q}}\mid
p,q\in SL_{2}\left( \mathbb{R}\right) \right\} \text{.}
\end{equation*}%
We set $\mathbb{R}\mathbf{1}=\mathbb{R}$, $\mathbb{C}=\mathbb{R}+\mathbf{i}%
\mathbb{R}$, $\mathbb{M}=\mathbb{C}+\mathbf{j}\mathbb{C}$ and write the
operators in $O_{o}\left( \mathbb{M}\right) $ in blocks with respect to this
decomposition. The ordered bases $\left\{ \mathbf{1},\mathbf{i}\right\} $
and $\left\{ \mathbf{j},\mathbf{ji}=-\mathbf{k}\right\} $ of $\mathbb{C}$
and $\mathbf{j}\mathbb{C}$ determine an isomorphism from $SO_{2}\cong S^{1}$
onto $SO\left( \mathbf{j}\mathbb{C},\mathbb{C}\right) $, assigning to $u\in
S^{1}$ the map $\mathbf{j}\mathbb{C}\ni \mathbf{j}z\mapsto uz$. The
juxtaposition of these bases provides the identification $O_{o}\left(
\mathbb{M}\right) \cong O_{o}\left( 2,2\right) $.

Now we consider the M\"{o}bius action $\ast :O_{o}\left( 2,2\right) \times
SO_{2}\rightarrow SO_{2}$, using the presentation in (\ref{action}),
identifying $u\in SO_{2}\equiv S^{1}$ with the set $\left\{ uz+\mathbf{j}%
z\mid z\in \mathbb{C}\right\} $, according to (\ref{Ugrasss}). Thus,%
\begin{equation}
L_{p}R_{\overline{q}}\ast u=v\left( p,q,u\right) \text{,}  \label{uvu}
\end{equation}%
where $v\left( p,q,u\right) $ is the only point of $S%
{{}^1}%
$ that satisfies that for all $z\in \mathbb{C}$ exists $w\left( z\right) \in
\mathbb{C}$ such that
\begin{equation*}
p\left( uz+\mathbf{j}z\right) \overline{q}=v\left( p,q,u\right) w\left(
z\right) +\mathbf{j}w\left( z\right) \text{.}
\end{equation*}

\begin{prp}
\label{O22}\emph{a)} The M\"{o}bius action $O_{o}\left( 2,2\right) \times
SO_{2}\rightarrow SO_{2}$ is not almost effective. More precisely, $\left(
L_{p}\circ R_{\overline{q}}\right) \ast u=L_{p}\ast u$ for all $p,q\in
SL_{2}\left( \mathbb{R}\right) $ and $u\in S^{1}$.

\smallskip

\emph{b)} The following diagram commutes\emph{:}%
\begin{equation*}
\begin{array}{ccc}
SL_{2}\left( \mathbb{R}\right) \times S^{1} & \overset{\cdot }{%
\longrightarrow } & S^{1} \\
\downarrow \left( F,\rho \right) &  & \downarrow \rho \\
O_{o}\left( 2,2\right) \times SO_{2} & \overset{\ast }{\longrightarrow } &
SO_{2}\text{,}%
\end{array}%
\end{equation*}%
where $F\left( p\right) =L_{p}$, $\rho \left( u\right) \left( z\right) =uz$
and $\cdot $ denotes the usual action by conformal transformations of the
plane preserving the circle, that is,
\begin{equation}
p\cdot u=\left( \alpha u+\overline{\beta}\right) \left( \beta u+\overline{%
\alpha}\right) ^{-1}  \label{pu}
\end{equation}%
if $p=\alpha +\mathbf{j}\beta $, with $\alpha ,\beta \in \mathbb{C}$, $%
\left\vert \alpha \right\vert ^{2}-\left\vert \beta \right\vert ^{2}=1$.
\end{prp}

	\begin{proofn}
a) Let $q=\alpha +\mathbf{j}\beta $ with $\alpha ,\beta \in \mathbb{C}$.
Given $u\in S^{1}$, we compute
\begin{equation*}
R_{\overline{\alpha}-\mathbf{j}\beta}\left( uz+\mathbf{j}z\right) =\left( uz+%
\mathbf{j}z\right) \left( \overline{\alpha}-\mathbf{j}\beta\right) =uz%
\overline{\alpha}-\overline{z}\beta+\mathbf{j}\left( z\overline{\alpha }-%
\overline{u}\,\overline{z}\beta\right)
\end{equation*}
(we have used that $w\mathbf{j}=\mathbf{j}\overline{w}$ for all $w\in
\mathbb{C}$). Since $uz\overline{\alpha}-\overline{z}\beta=u\left( z%
\overline{\alpha}-\overline{u}\,\overline{z}\beta\right) $ for all $z\in
\mathbb{C}$, we have by (\ref{uvu}) that $v\left( 1,q,u\right) =u$, and so, $%
R_{\overline{q}}\ast u=u$.

\medskip

b) Now let $p=\alpha+\mathbf{j}\beta\in SL_{2}\left( \mathbb{R}\right) $
with $\alpha,\beta\in\mathbb{C}$. We compute%
\begin{equation*}
L_{\alpha+\mathbf{j}\beta}\left( uz+\mathbf{j}z\right) =\left( \alpha+%
\mathbf{j}\beta\right) (uz+\mathbf{j}z)=\left( \alpha u+\overline {\beta}%
\right) z+\mathbf{j}\left( \beta u+\overline{\alpha}\right) z\text{.}
\end{equation*}
Calling $w=\left( \beta u+\overline{\alpha}\right) z$, we have $z=\left(
\beta u+\overline{\alpha}\right) ^{-1}w$, and so $v\left( p,1,u\right) $ as
in (\ref{uvu}) equals the right hand side of (\ref{pu}), as desired.
\qed
\end{proofn}

\subsection{Relationship with conformal motions of the three sphere}


\begin{lemn}
Let $f:sp\left( 1,1\right) \rightarrow o\left( 1,4\right) $ be defined by
\begin{equation*}
f\left(
\begin{array}{cc}
\xi & \overline{\alpha } \\
\alpha & \eta%
\end{array}%
\right) =\left(
\begin{array}{cc}
0 & 2\alpha ^{T} \\
2\alpha & L_{\eta }-R_{\xi }%
\end{array}%
\right) ,
\end{equation*}%
where $\xi ,\eta \in \operatorname{Im}(\mathbb{H})$ and $\alpha \in \mathbb{H}$.
Then $f$ is a Lie algebra isomorphism which induces a morphism $F:Sp\left(
1,1\right) \rightarrow O_{o}\left( 1,4\right) $.
\end{lemn}

\begin{proofn}
Standard arguments show that $f$ is a morphism, using that%
\begin{equation*}
4\left( \alpha \beta ^{T}-\beta \alpha ^{T}\right) =L_{\alpha \overline{%
\beta }-\beta \overline{\alpha }}-R_{\overline{\alpha }\beta -\overline{%
\beta }\alpha }
\end{equation*}%
holds for all $\alpha ,\beta \in \mathbb{H}$ (in the left hand side,
quaternions are thought of as column vectors, under the canonical
identification of $\mathbb{H}$ with $\mathbb{R}^{4}$). Now we sketch how to
verify this identity. Since each side defines an skew symmetric $\mathbb{R}$%
-bilinear map
\begin{equation*}
\left( \alpha ,\beta \right) \in \mathbb{H}\times \mathbb{H}\rightarrow
\left\{ \text{skew-symmetric }\mathbb{R}\text{-linear maps of }\mathbb{H}%
\right\} \text{,}
\end{equation*}%
it suffices to consider\ only the case when $\left\{ \alpha ,\beta \right\} $
is orthonormal. By linearity, it is enough to evaluate both sides at $%
q=\alpha $, $q=\beta $ and $q\in \left\{ \alpha ,\beta \right\} ^{\bot }$.

The last assertion follows, since $Sp\left( 1,1\right) $ is simply connected
(by the Cartan decomposition, it is diffeomorphic to $\left( S^{3}\times
S^{3}\right) \times \mathbb{R}^{6}$).
\qed
\end{proofn}

\begin{prp}
\label{conmuSp11}The following diagram is commutative\emph{:}%
\begin{equation*}
\begin{array}{ccc}
Sp\left( 1,1\right) \times Sp_{1} & \overset{\ast }{\longrightarrow } &
Sp_{1} \\
\downarrow \left( F,\rho \right) &  & \downarrow \rho \\
O_{o}\left( 1,4\right) \times S^{3} & \overset{\ast _{c}}{\longrightarrow }
& S^{3}\text{,}%
\end{array}%
\end{equation*}%
where $\rho \left( u\right) =\overline{u}$, $\ast $ denotes the M\"{o}bius
action and $\ast _{c}$ the conformal action.
\end{prp}

\begin{proofn}
First, we evaluate $F$ at the elements of the maximal compact subgroup $%
K=Sp_{1}\times Sp_{1}$ of $Sp\left( 1,1\right) $. For all $p,q,u\in S^{3}$,
\begin{equation*}
F\left( \operatorname{diag}\left( p,q\right) \right) \ast _{c}\rho \left( u\right)
=L_{q}R_{\overline{p}}\ast _{c}\overline{u}=q\,\overline{u}\,\overline{p}%
=\rho \left( \operatorname{diag}\left( p,q\right) \ast u\right)
\end{equation*}%
holds, putting $q=e^{t\eta }$ for $\eta \in $ $\operatorname{Im}(\mathbb{H})$, since
$\exp \left( tL_{\eta }\right) =L_{e^{t\eta }}$ (and similarly for $R_{%
\overline{p}}$).

Given that $F$ is a morphism, by the $KAK$ decomposition of $Sp\left(
1,1\right) $, where $A=\exp \left( \mathbb{R}X\right) $ with $X=\left(
\begin{array}{cc}
0 & 1 \\
1 & 0%
\end{array}%
\right) $, it only remains to check that
\begin{equation*}
F\left( \exp tX\right) \ast _{c}\rho \left( u\right) =\rho \left( \exp
\left( tX\right) \ast u\right)
\end{equation*}%
holds for all $t\in \mathbb{R}$. We call
\begin{equation*}
H_{t}=\exp \left( tX\right) =\left(
\begin{array}{cc}
\cosh t & \sinh t \\
\sinh t & \cosh t%
\end{array}%
\right) \in O_{o}\left( 1,1\right) \subset Sp\left( 1,1\right) \text{.}
\end{equation*}%
By definition of the M\"{o}bius action we have%
\begin{equation}
\rho \left( \exp \left( tX\right) \ast u\right) =\rho \left( H_{t}\ast
u\right) =\left( \overline{u}\sinh t+\cosh t\right) ^{-1}\left( \overline{u}%
\cosh t+\sinh t\right) \text{.}  \label{sp11sb}
\end{equation}

On the other hand,
\begin{equation*}
F\left( \exp tX\right) =\exp \left( tf\left( X\right) \right) =\exp \left( t%
\operatorname{diag}\left( 2X,0_{3}\right) \right) =\operatorname{diag}\left(
H_{2t},I_{3}\right) \text{.}
\end{equation*}%
Now, putting $u=a+w\in S^{3}$, with $a\in \mathbb{R}$, $w\in \operatorname{Im}\left(
\mathbb{H}\right) $, we compute
\begin{equation*}
\left(
\begin{array}{ccc}
\cosh 2t & \sinh 2t & 0 \\
\sinh 2t & \cosh 2t & 0 \\
0 & 0 & I_{3}%
\end{array}%
\right) \left(
\begin{array}{c}
1 \\
a \\
-w%
\end{array}%
\right) =\left(
\begin{array}{c}
\cosh \left( 2t\right) +\sinh \left( 2t\right) a \\
\sinh \left( 2t\right) +\cosh \left( 2t\right) a \\
-w%
\end{array}%
\right) \text{.}
\end{equation*}%
So, by (\ref{conformal}),
\begin{equation}
F\left( \exp tX\right) \ast _{c}\overline{u}=\operatorname{diag}\left(
H_{2t},I_{3}\right) \ast _{c}\overline{u}=\frac{\sinh \left( 2t\right)
+a\cosh \left( 2t\right) -w}{\cosh (2t)+a\sinh (2t)}\text{.}  \label{Sp11con}
\end{equation}

Using that $\left\vert u\right\vert ^{2}=a^{2}+\left\vert w\right\vert
^{2}=1 $ and $1=u\overline{u}=\left( a+w\right) \left( a-w\right)
=a^{2}-w^{2}$, a lengthy but straightforward computation yields that
expressions (\ref{sp11sb}) and (\ref{Sp11con}) coincide, as desired.
\qed\end{proofn}

\subsection{Relationship with projective motions of the three sphere}


Identifying $\mathbb{R}^{4}=\mathbb{H}$ and $\mathbb{R}^{3}=\operatorname{Im}(%
\mathbb{H})$, the orthogonal groups for dimensions 3 and 4 can be described
as
\begin{equation*}
SO_{3}=\left\{ L_{p}\circ R_{\overline{p}}\mid p\in S^{3}\right\} \text{ \ \
\ \ \ \ and \ \ \ \ \ \ \ }SO_{4}=\left\{ L_{p}\circ R_{\overline{q}}\mid
p,q\in S^{3}\right\} \text{.}
\end{equation*}

Next we recall the well known double covering morphism from $SL_{4}(\mathbb{R%
})$ to $O_{o}\left( 3,3\right) $. Let $o_{4}=\left\{ \text{skew symmetric }%
4\times 4\text{ real matrices}\right\} $. There is an isomorphism $\omega
:o_{4}\rightarrow \Lambda ^{2}\left( \mathbb{R}^{4}\right) ^{\ast }$, $%
Z\mapsto \omega _{Z}$, where $\omega _{Z}\left( x,y\right) =x^{T}Zy$. We
consider on $o_{4}$ the inner product $g$ of signature $\left( 3,3\right) $
induced via $\omega $ from the inner product $h$ on $\Lambda ^{2}\left(
\mathbb{R}^{4}\right) ^{\ast }$ given by $h\left( \sigma ,\tau \right) {%
\mathrm{det}}=\frac{1}{2}\,\sigma \wedge \tau $.

We also have the decomposition
\begin{equation*}
o_{4}=\mathcal{L}\oplus \mathcal{R}\cong \operatorname{Im}(\mathbb{H})\times (%
\mathbb{H})\text{,}
\end{equation*}%
where $\mathcal{L}=\left\{ L_{\xi }\mid \xi \in \operatorname{Im}(\mathbb{H}%
)\right\} $ and $\mathcal{R}=\left\{ R_{\eta }\mid \eta \in \operatorname{Im}(%
\mathbb{H})\right\} $. The set $\mathcal{B}=\left\{
L_{i},L_{j},L_{k}\right\} \cup \left\{ R_{i},R_{j},R_{k}\right\} $ is a $%
\left( 3,3\right) $-orthonormal basis of $o_{4}$ with respect to $g$, where
the first elements are space-like and the last ones are time-like. Now, the
group $SL_{4}\left( \mathbb{R}\right) $ acts on $\Lambda ^{2}\left( \mathbb{R%
}^{4}\right) ^{\ast }$ by linear isometries with respect to $h$ via $%
A.\sigma =\left( A^{T}\right) ^{\ast }\sigma $. This action, through the
isomorphism $\omega ,$ induces the morphism
\begin{equation*}
\Phi :SL_{4}(\mathbb{R})\rightarrow O_{o}\left( o_{4},g\right) \text{,\ \ \
\ \ \ \ \ \ \ \ }\Phi \left( A\right) \left( Z\right) =AZA^{T}\text{.}
\end{equation*}%
By means of the basis $\mathcal{B}$, we identify $o_{4}=\mathbb{R}^{3,3}$
and so $\Phi $ induces a morphism $F:SL_{4}(\mathbb{R})\rightarrow
O_{o}\left( 3,3\right) $. We also define $\rho :S^{3}\rightarrow SO_{3}$ by $%
\rho \left( q\right) =L_{q}R_{\overline{q}}$.

\begin{prp}
\label{conmuSL4}The following diagram is commutative\emph{:}%
\begin{equation}
\begin{array}{ccc}
SL_{4}(\mathbb{R})\times S^{3} & \overset{\ast _{p}}{\longrightarrow } &
S^{3} \\
\downarrow \left( F,\rho \right) &  & \downarrow \rho \\
O_{o}\left( 3,3\right) \times SO_{3} & \overset{\ast }{\longrightarrow } &
SO_{3}\text{,}%
\end{array}
\label{diagramaSL4}
\end{equation}%
where $\ast $ denotes the M\"{o}bius action and $\ast _{p}$ the projective
action.
\end{prp}

\begin{proofn}
First we verify the commutativity of the diagram for the elements of $%
SO_{4}\subset SL_{4}\left( \mathbb{R}\right) $. We have that $\left.
F\right\vert _{SO_{4}}:SO_{4}\rightarrow SO_{3}\times SO_{3}\subset
O_{o}\left( 3,3\right) $ (canonical immersion) is given by
\begin{equation*}
F\left( L_{p}\circ R_{\overline{q}}\right) =\left( L_{p}\circ R_{\overline{p}%
},L_{q}\circ R_{\overline{q}}\right) \text{.}
\end{equation*}%
Indeed, calling $L_{p}\circ R_{\overline{q}}=A$, $A\ast _{p}u=pu\overline{q}$
holds for all $u\in S^{3}$, since $\left\vert pu\overline{q}\right\vert =1$.
Now, $F\left( A\right) \ast \rho \left( u\right) =\rho \left( A\ast
_{p}u\right) $ if and only if
\begin{equation*}
L_{p}R_{\overline{p}}L_{u}R_{\overline{u}}\left( L_{q}R_{\overline{q}%
}\right) ^{-1}=L_{pu\overline{q}}\left( R_{pu\overline{q}}\right) ^{-1}
\end{equation*}%
for all $p,q,u\in S^{3}$, and this can be easily checked.

Since $F$ is a morphism, by the $KAK$ decomposition of $SL_{4}\left( \mathbb{%
R}\right) $, it remains to check that the diagram commutes for any diagonal
matrix $E$ in $SL_{4}\left( \mathbb{R}\right) $. Any such $E$ has the form $%
E=\exp \left( x_{1}D+x_{2}D_{2}+x_{3}D_{3}\right) $, where 
\begin{equation*}
D=\operatorname{diag}%
\left( 1,1,-1,-1\right), \ \ D_{2}=\operatorname{diag}\left( 1,-1,1,-1\right) 
\text{\ \ \ \ and \ \ \ }
D_{3}=\operatorname{diag}\left( 1,-1,-1,1\right). 
\end{equation*}
Now, there exist $k_{2},k_{3}\in
SO_{4}$ such that $D_{l}=k_{l}Dk_{l}^{-1}$ for $l=2,3.$ Hence,
\begin{equation*}
E=\exp \left( x_{1}D\right) k_{2}\exp \left( x_{2}D\right)
k_{2}^{-1}k_{3}\exp \left( x_{1}D\right) k_{3}^{-1}.
\end{equation*}%
Therefore, it suffices to check the commutativity of the diagram for
\begin{equation*}
A_{t}:=\exp \left( tD\right)= \operatorname{diag}\left(
e^{t},e^{t},e^{-t},e^{-t}\right).
\end{equation*}

Let $B_{t}$ be the matrix of the transformation $\Phi \left( A_{t}\right)
\left( L_{i}\right) =A_{t}L_{i}A_{t}^{T}\in o\left( \mathbb{H}\right) $ with
respect to the basis $\left\{ 1,i,j,k\right\} $. We compute%
\begin{eqnarray*}
B_{t} &=&\left(
\begin{array}{cc}
e^{t}I_{2} & 0_{2} \\
0_{2} & e^{-t}I_{2}%
\end{array}%
\right) \left(
\begin{array}{cc}
J & 0 \\
0 & J%
\end{array}%
\right) \left(
\begin{array}{cc}
e^{t}I_{2} & 0_{2} \\
0_{2} & e^{-t}I_{2}%
\end{array}%
\right) =\left(
\begin{array}{cc}
e^{2t}J & 0_{2} \\
0_{2} & e^{-2t}J%
\end{array}%
\right) \\
&=&\cosh \left( 2t\right) \left(
\begin{array}{cc}
J & 0_{2} \\
0_{2} & J%
\end{array}%
\right) +\sinh \left( 2t\right) \left(
\begin{array}{cc}
J & 0_{2} \\
0_{2} & -J%
\end{array}%
\right)
\end{eqnarray*}%
for all $t\in \mathbb{R}$, where $J=\left(
\begin{array}{cc}
0 & -1 \\
1 & 0%
\end{array}%
\right) $. Similar computations yield%
\begin{eqnarray*}
A_{t}L_{ai+bj+ck}A_{t}^{T} &=&\cosh \left( 2t\right) L_{ai}+\sinh \left(
2t\right) R_{ai}+L_{bj+ck} \\
A_{t}R_{ai+bj+ck}A_{t}^{T} &=&\sinh \left( 2t\right) L_{ai}+\cosh \left(
2t\right) R_{ai}+R_{bj+ck}
\end{eqnarray*}%
for all $a,b,c,t\in \mathbb{R}$. Consequently, calling
\begin{equation*}
C_{t}=\operatorname{diag}\left( \cosh 2t,I_{2}\right) \text{\ \ \ \ \ and\ \ \ \ \ }%
S_{t}=\operatorname{diag}\left( \sinh 2t,0_{2}\right) \text{,}
\end{equation*}%
both in $\mathbb{R}^{3\times 3}$, the matrix of the operator $\Phi \left(
A_{t}\right) $ with respect to the basis $\mathcal{B}$, which equals the
matrix of $F\left( A_{t}\right) $, is
\begin{equation*}
\left(
\begin{array}{cc}
C_{t} & S_{t} \\
S_{t} & C_{t}%
\end{array}%
\right) \in O_{o}\left( 3,3\right) \text{.}
\end{equation*}%
Now, for $w\neq 0$ in $\mathbb{H}$, we call $I_{w}$ the matrix of the
operator $L_{w}\left( R_{w}\right) ^{-1}$ on $\operatorname{Im}(\mathbb{H})$ with
respect to the basis $\left\{ i,j,k\right\} $. If $w=a+bi+cj+dk$, we compute%
\begin{equation*}
I_{w}=\frac{1}{\left\vert w\right\vert ^{2}}\left(
\begin{array}{ccc}
a^{2}+b^{2}-d^{2}-c^{2} & -2ad+2bc & 2bd+2ac \\
2ad+2bc & a^{2}-b^{2}+c^{2}-d^{2} & -2ab+2cd \\
-2ac+2bd & 2ab+2cd & a^{2}-b^{2}+d^{2}-c^{2}%
\end{array}%
\right) \text{.}
\end{equation*}%
The commutativity of the diagram for $A_{t}$ amounts to the validity of%
\begin{equation*}
\left( C_{t}I_{u}+S_{t}\right) \left( S_{t}I_{u}+C_{t}\right)
^{-1}=L_{\left( A_{t}u\right) /\left\vert A_{t}u\right\vert }\left(
R_{\left( A_{t}u\right) /\left\vert A_{t}u\right\vert }\right)
^{-1}=L_{A_{t}u}\left( R_{A_{t}u}\right) ^{-1}=I_{A_{t}u}
\end{equation*}%
for all $u\in S^{3}$ and all $t\in \mathbb{R}$, or equivalently, $%
C_{t}I_{u}+S_{t}=I_{A_{t}u}\left( S_{t}I_{u}+C_{t}\right) $. This can be
checked with a lengthy but straightforward computation.
\qed
\end{proofn}

\subsection{Dynamical implications for M\"{o}bius motions in low dimensions}

Let $Sp\left( 1,1\right) $ and $O_{o}\left( 3,3\right) $ be endowed with the
kinetic energy metrics. By Proposition \ref{CompactoMaximal}, the induced
Riemannian metrics on the maximal compact subgroups $Sp_{1}\times Sp_{1}$
and $SO_{3}\times SO_{3}$ are bi-invariant. Theorem \ref{XceroceroY} asserts
that each factor is totally geodesic. The results in this section give more
information on the matter. At the same time, they provide an alternative
statement (perhaps more illuminating) for Theorem 3 (b) in \cite{LSW}.

\begin{cor}
\label{coro} Let $Sp(1,1)$ and $O_{o}\left( 3,3\right) $ be endowed with the
kinetic energy metrics. Then:

\smallskip

\emph{a)} The subgroup $Sp_{1}\times Sp_{1}$ is totally geodesic in $%
Sp\left( 1,1\right) $.

\smallskip

\emph{b)} A geodesic in $SO_{3}\times SO_{3}$ is a geodesic in $O_{o}\left(
3,3\right) $ if and only if it remains in one of the factors.
\end{cor}

\begin{proofn}
By Propositions \ref{conmuSp11} and \ref{conmuSL4}, the M\"{o}bius actions
of $Sp\left( 1,1\right) $ on $Sp_{1}$ and of $O_{o}\left( 3,3\right) $ on $%
SO_{3}$ are equivalent, up to double coverings, to the conformal action of $%
O_{o}\left( 1,4\right) $ on $S^{3}$ and the projective action of $%
SL_{4}\left( \mathbb{R}\right) $ on $S^{3}$, respectively. The morphisms
preserve the corresponding kinetic energy metrics.

By Theorem 2 (b) in \cite{Sal}, $SO_{4}$ is totally geodesic in $O_{o}\left(
1,4\right) $ endowed with the kinetic energy metric associated with the
conformal action. Since $Sp_{1}\times Sp_{1}$ corresponds to $SO_{4}$ under
the corresponding morphism, part (a) of the corollary follows.

Also, Theorem 3 (b) in \cite{LSW} asserts that for $Z\in o_{4}$, the
geodesic $\sigma (t)=\exp (tZ)$ of $SO_{4}$ through the identity is also a
geodesic of $SL_{4}\left( \mathbb{R}\right) $ if and only if $Z=\lambda X$
for some $\lambda \in \mathbb{R}$ and $X\in o_{4}$ with $X^{2}=-I_{4}$. Now,
$Z=L_{\xi }+R_{\eta }$ for some $\xi ,\eta \in \operatorname{Im}(\mathbb{H})$. We
compute%
\begin{equation*}
-\lambda ^{2}I_{4}=Z^{2}=\left( L_{\xi }+R_{\eta }\right) ^{2}=L_{\xi
}^{2}+R_{\eta }^{2}+2L_{\xi }R_{\eta }=-\left( \left\vert \xi \right\vert
^{2}+\left\vert \eta \right\vert ^{2}\right) I_{4}+2L_{\xi }R_{\eta }
\end{equation*}%
($L_{\xi }$ and $R_{\eta }$ commute). This implies that $L_{\xi }R_{\eta
}=\mu I_{4}$ for some $\mu \in \mathbb{R}$. If $\eta \neq 0$, we have $%
L_{\xi }=R_{\mu \bar{\eta}/\left\vert \eta \right\vert ^{2}}$. Now, $p,q\in
\mathbb{H}$ satisfy $L_{p}=R_{q}$ only if $p,q$ are in the center of $%
\mathbb{H}$, that is, $p,q\in \mathbb{R}$. Since $\xi \in \operatorname{Im}(\mathbb{H%
})$, then $\xi =0$, as desired.
\qed
\end{proofn}

\medskip

\medskip
\noindent {\bf Acknowledgements.} \textsl{This work was supported by Consejo Nacional de Investigaciones
Cient\'ificas y T\'ecnicas and Secretar\'{\i}as de Ciencia y T\'ecnica of Universidad Nacional de C\'ordoba  and  Universidad Nacional de Rosario.
We would like to thank Carlos Olmos for
useful comments and Simon Chiossi for a valuable reference.}


%
%
%
%
%


\begin{thebibliography}{99}
	\bibitem{ArnoldArt} Arnold, V.\,I., Sur la g\'eom\'etrie
diff\'erentielle des groupes de Lie  dimension infinie et ses applications
a l'hydrodynamique des fluids parfaits, \textit{Ann. Inst. Grenoble}, 1966, vol.\,16,
pp.\,316--361.

\bibitem{BetterThanNice} Bates, L. and Gibson, P., A geometry where everything
is better than nice, \textit{Proc. Amer. Math. Soc.}, 2017, vol.\,145, pp.\,461--465.

\bibitem{BBM14} Bauer, M., Bruveris, M. and Michor P.\,W.,
Overview of the geometries of shape spaces and diffeomorphism
groups, \textit{J. Math. Imaging Vis.}, 2014, vol.\,50, pp.\,60--97.

\bibitem{BauerModinSemi} Bauer, M. and Modin, K., Semi-invariant Riemannian
metrics in hydrodynamics, \textit{Calc. Var.}, 2020, vol.\,59, paper
No. 65, 25 pp.

\bibitem{Olmos} Berndt, J., Console, S. and Olmos, C., \textit{Submanifolds and
	Holonomy}, Monographs and Research Notes in Mathematics, Boca
Raton, FL: CRC Press, 2016.

\bibitem{BruverisVialardJEMS} Bruveris, M. and Vialard, F.-X., On
completeness of groups of diffeomorphisms, \textit{J. Eur. Math. Soc.}, 2017, vol.\,19,
pp.\,1507--1544.

\bibitem{Burov} Burov, A.\,A. and Chevallier, D.\,P., Dynamics of affinely
deformable bodies from the standpoint of theoretical mechanics and
differential geometry, \textit{Rep. Math. Phys.}, 2008, vol.\,62, pp.\,283--321.

\bibitem{CapSlovak} \v{C}ap, A. and Slov\'ak, J., \textit{%
	Parabolic Geometries. I. Background and general theory}, Mathematical
Surveys and Monographs 154, Providence, RI: American Mathematical Society,
2009.

\bibitem{EbinSympl} Ebin, D.\,G., Geodesics on the symplectomorphism
group, \textit{Geom. Funct. Anal.}, 2012, vol.\,22, pp.\,202--212.

\bibitem{EbinContact} Ebin, D.\,G. and Preston, S.\,C., Riemannian geometry of
the contactomorphism group, \textit{Arnold Math. J.}, 2015, vol.\,1, pp.\,5--36.

\bibitem{EbinMarsden} Ebin, D.\,G. and Marsden, J., Groups of diffeomorphisms and
the motion of an incompressible fluid, \textit{Ann. of Math.}\,(2), 1970, vol.\,9, pp.\,102--163.

\bibitem{ES} Emmanuele, D. and Salvai, M., Force free M\"{o}bius motions of the
circle, \textit{J. Geom. Symmetry Phys},  2012, vol.\,27, pp.\,59--65.

\bibitem{ESarXiv} Emmanuele, D. and Salvai, M., A dynamical presentation of the
better than nice metric on the disc, arXiv:1605.03863 [math.DG] (postprint 2016).

\bibitem{Viz02} Haller, S., Teichmann, J. and Vizman, C., Totally geodesic
subgroups of diffeomorphisms, \textit{J. Geom. Phys.}, 2002, vol.\,42, pp.\,342--354.

\bibitem{Har} Harvey, F.\,R., \textit{Spinors and Calibrations}, Perspectives in
Mathematics 9, Boston, MA: Academic Press,  1990.

\bibitem {HelG}Helgason, S., \textit{Differential Geometry, Lie Groups and
	Symmetric Spaces}, Graduate Studies in Mathematics 34, Providence, RI: American
Mathematical Society,  2001.


\bibitem{Khesin} Khesin, B., Misio\l ek, G. and Modin, K. \textsl{Geometric
	hydrodynamics and infinite-dimensional Newton's equations}. arXiv:2001.01143
[math.SG] (2020), to appear in Bulletin of the Amer. Math. Soc.

\bibitem{Kh2} Khesin, B., Dynamics of symplectic fluids and point
vortices, \textit{Geom. Funct. Anal.},  2012, vol.\,22, pp.\,1444--1459.

\bibitem{LSW} Lazarte, M., Salvai, M. and Will, A., Force free projective
motions of the sphere, \textit{J. Geom. Phys.},  2007, vol.\,57, pp.\,2431--2436.

\bibitem{Me} Meckes, E., \textit{The Random Matrix Theory of the Classical
	Compact Groups}, Cambridge Tracts in Mathematics 218, Cambridge: Cambridge
University Press,  2019.

\bibitem{ModinPerl} Modin, K., Perlmutter, M., Marsland, S. and McLachlan, R.\,I.,
On Euler-Arnold equations and totally geodesic subgroups, \textit{J. Geom.
	Phys.}, 2011, vol.\,61, pp.\,1446--1461.

\bibitem{Nagano} Nagano, T., Transformation groups on compact symmetric spaces,
 \textit{Trans. Amer. Math. Soc.}, 1965, vol.\,118, pp.\,428--453.

\bibitem{ONeill} O'Neill, B., \textit{Semi-Riemannian Geometry (with
	applications to Relativity)}, New York: Academic Press, 1983.

\bibitem{OK87} Ovsienko, V.\,Y. and Khesin, B., Korteweg--de Vries
superequations as an Euler equation, \textit{Funct. Anal. Appl.},
1987, vol.\,21, pp.\,329--331.

\bibitem{Sal} Salvai, M., Force free conformal motions of the sphere,  \textit{Diff.
	Geom. Appl.},  2002, vol.\,16, pp.\,285--292.

\bibitem{STohoku} Salvai, M., Circles in self-dual symmetric $R$%
	-spaces, \textit{Tohoku Math. J.} (2), 2021, vol.\,73,  pp.\,257-275.

\bibitem{Sharpe} Sharpe, R., \textit{Differential Geometry: Cartan's
	generalisation of Klein's Erlangen program}, Graduate Texts in Mathematics
166, New York: Springer-Verlag,  1997.

\bibitem{Sideris} Sideris, T.\,C., Global existence and asymptotic behavior of affine motion of 3D ideal fluids surrounded by vacuum, \textit{Arch. Ration. Mech. Anal.},  2017, vol.\,225, pp.\,141--176.

\bibitem{Si} Simon, B., \textit{Representations of Finite and Compact Groups},
Graduate Studies in Mathematics 10, Providence, RI: American Mathematical Society,
1996.

\bibitem{Ta} Takeuchi, M., Basic transformations of symmetric $R$-spaces, \textit{Osaka J. Math.},  1988, vol.\,25, pp.\,259--297.

\bibitem{Th} Thorbergsson, G., Classical symmetric $R$-spaces, \textit{Rend. Semin.
	Mat. Univ. Politec. Torino}, 2016, vol.\,74, pp.\,329--354.

\bibitem{Viz08} Vizman, C., Geodesic equations on diffeomorphism
groups, \textit{SIGMA Symmetry Integrability Geom. Methods Appl.}, 2008, vol.\,4, paper
No. 030, 22 pp.

\bibitem{Zh} Zhang, F., Quaternions and matrices of quaternions, \textit{Linear
	Algebra Appl.}, 1997, vol.\,251, pp.\,21--57.
\end{thebibliography}
\end{document}